\documentclass[11pt]{article}

\usepackage{amsmath}
\usepackage{graphicx}
\usepackage{amssymb}
\usepackage{epstopdf}
\usepackage{epsfig}
\usepackage{subfigure}
\usepackage{amsmath, amssymb, amsfonts}
\usepackage{a4, a4wide}
\usepackage{color}
\usepackage{verbatim}
\usepackage{array}
\usepackage{hyperref}
\usepackage{tikz}
\usetikzlibrary{matrix}
\def\ie{{\rm i.e.,\/}\ }
\def\etc{{\rm etc.\/}\ }

     \newcommand{\und}[1]{\underline{#1}} 
\newcommand{\nco}{\newcommand}
\nco{\one}{\ensuremath{\,\,\mathrm{l}\!\!\!1}} 
\nco{\NN}{\mathbb{N}}
\nco{\ZZ}{\mathbb{Z}}
\nco{\QQ}{\mathbb{Q}}
\nco{\RR}{\mathbb{R}}
\nco{\CC}{\mathbb{C}}
\nco{\HH}{\mathbb{H}}
\nco{\OO}{\mathbb{O}}
\nco{\red}{\color{red}}
\nco{\redend}{\normalcolor}
\nco{\0}{\underline{0}}
\nco{\1}{\underline{1}}
\nco{\2}{\underline{2}}
\nco{\3}{\underline{3}}
\DeclareMathOperator{\sfac}{sf}
\newtheorem{theorem}{Theorem}

\newtheorem{definition}{Definition}

\makeatletter
\newcommand{\pmatrixcmd}[1]{}%
\DeclareRobustCommand{\pmatrixcmd}[1]{\left(\env@matrix#1\endmatrix\right)}
\makeatother

\begin{document}

\noindent
 
 \vspace{0.2cm} 
 
\begin{center}
\begin{Large}
    {Quantum McKay correspondence and global dimensions for fusion and module-categories associated with Lie groups}
\end{Large}
\end{center}


\begin{center}
\renewcommand{\thefootnote}{\arabic{footnote}}
{R. Coquereaux} \footnotemark[1]${}^,$ \footnotemark[2]
\\ \renewcommand{\thefootnote}{\arabic{footnote}} 
\end{center}

 \vspace{0.2cm} 
 

 \vspace{0.3cm} 

\abstract{Global dimensions for fusion categories ${\mathcal A}_k(G)$ defined by a pair $(G,k)$, where $G$ is a Lie group and $k$ a positive integer, are expressed in terms of Lie quantum superfactorial functions. The global dimension is defined as the square sum of quantum dimensions of simple objects, for the category of integrable modules over an affine Lie algebra at some level. The same quantities can also be defined from the theory of quantum groups at roots of unity or from conformal field theory WZW models.
Similar results  are also presented for those associated module-categories that can be obtained via conformal embeddings (they are ``quantum subgroups'' of a particular kind).
As a side result, we express the classical (or quantum) Weyl denominator of simple Lie groups in terms of classical (or quantum) factorials calculated for the exponents of the group.
Some calculations use the correspondence existing between periodic quivers for simply-laced Lie groups and fusion rules for  module-categories associated with ${\mathcal A}_k(SU(2))$.}

\vspace{0.3cm}


\vspace{0.2cm}

\noindent{{Classification}}: 81R50; 81T40; 18D10; 33E99. 

\addtocounter{footnote}{0}
\footnotetext[1]{{\scriptsize{\it IMPA, Instituto Nacional de Matem\'atica Pura e Aplicada, Rio de Janeiro, Brasil}.}}
\footnotetext[2]{{\scriptsize{\it UMI CNRS-IMPA (2924), on leave from Aix Marseille Universit\'e, CNRS, CPT, UMR 7332, 13288 Marseille, France}.}}
\vspace{0.4cm}
\unitlength = 1mm
\section{Introduction}

Let $k$ be a positive integer and $G$ a simple, compact Lie group and $Lie(G)$ its complex Lie algebra.
We call ${\mathcal A}_k(G)$ the category of integrable modules of the Kac-Moody algebra associated with 
$Lie(G)$ at level $k$, see e.g. \cite{Kac:book}. It is semi-simple, monoidal (use the fusion product), and modular.
It is also equivalent \cite{Finkelberg, KazhdanLusztig} to a category constructed in terms of  representations of the quantum group $G_q$ at  root of unity $q =  exp(\frac{i \pi}{N^\vee+k})$, where $N^\vee$ is the dual Coxeter number of $G$  (take the quotient of the category of tilting modules by the additive subcategory generated by indecomposable modules of zero quantum dimension).
These categories play a key role in the Wess -Ð Zumino -Ð Witten models of conformal field theory.

One can associate a quantum dimension  (or Perron-Frobenius dimension) $\mu$ to the objects of ${\mathcal A}_k(G)$,  
and since there are finitely many simple objects $n$, we can sum the square of their quantum dimensions $\mu(n)$ to define a global dimension $\vert {\mathcal A}_k(G) \vert = \sum_n \mu(n)^2$. This quantity generalizes the notion of order for a finite group.
There is  an explicit formula for this number, 
it involves a quantum version of the superfactorial where the product is done over q-integers associated with the exponents of the chosen Lie group. This function, that we call quantum Lie superfactorial of type $G$,  is therefore a kind of generalization of the quantum Barnes function (for the latter, just take $G=SU(N)$). 

One can also find an explicit formula for the global dimension $\vert {\mathcal E} \vert$ of a large class of module-categories ${\mathcal E}$ associated with ${\mathcal A}_k(G)$, 
namely those modules that can be obtained via conformal embedding  of $G$, at level $k$ in a larger Lie group $J$, at level 1 (see the text for details).
Notice that, using another language,  a module-category  ${\mathcal E}$ is a ``quantum subgroup'' of the quantum group $G_q$, for $q$ a root of unity.  

Both formulae will be recalled. They were presented  at a conference held in Cordoba (Argentina), in 2010. 
The article  \cite{CordobaRC},  from the proceedings, gives proofs and  examples, but does not detail the calculations leading to our definitions for the quantum Lie superfactorial of type $G$, in particular when the latter is exceptional.
Our proof of the formula for $\vert {\mathcal A}_k(G) \vert$ uses several known properties that will be recalled in the text but also an expression of the Weyl denominator (classical or quantum)
written in terms of product of factorials (classical or quantum) associated with the exponents of the Lie group $G$. 
As we did not find this expression in the literature, we also present proofs of the corresponding closed formulae obtained for the Weyl denominator, for the classical series of simple Lie groups as well as for the exceptional ones, simply-laced or not.  Those results can themselves be obtained in several ways. 
Remembering that the classical Weyl denominator is defined as the product, over all positive roots, of the scalar products between the latters and the Weyl vector, 
one possibility to obtain this family of scalar products is to perform a brute-force calculation using only standard Lie group theoretical properties. 
This approach will be followed in the Appendix for the classical series $A_r$, $B_r$, $C_r$ and $D_r$, for arbitrary values of the rank $r$. 
However, another method, less elementary but more in the spirit of the paper, relies on the quantum McKay correspondence relating ordinary ADE Dynkin diagrams $G$ and ``quantum subgroups of type $SU(2)$'', \ie module-categories over ${\mathcal A}_k(SU(2))$. This correspondence  has been known for many years and goes back to \cite{CIZ}, but here we use it as a tool, in the sense that  
we exploit a little-known relation previously described in \cite{Dorey:CoxeterElement} and  \cite{Ocneanu:Bariloche} between periodic quivers of roots for $G$ (those quivers are described for instance in \cite{KirillovThind}) and an underlying quantum $SU(2)$ theory.
Using this relation, and for $G$ of type $ADE$,  one can determine the scalar products between roots and the Weyl vector from the structure coefficients (fusion numbers) describing modules over the Grothendieck ring associated with the category ${\mathcal A}_k(SU(2))$. The non simply-laced cases are then handled by an appropriate modification of the same method.
This approach was briefly mentioned  and used in \cite{CordobaRC}, but no description of the corresponding methods is  available.  For this reason, and apart from a determination of the Weyl denominator leading, in section \ref{sec:periodicquiver}, to various classical or quantum Lie superfactorial functions and therefore to the closed formulae displayed in the first part,   the discussion to be found in section \ref{McKayCorrespondence} can be of independent interest.

Section \ref{sec: globaldim} summarizes the main results concerning  the global dimension of the fusion category ${\mathcal A}_k(G)$ and (a particular family of) associated modules-categories ${\mathcal E}$.
In particular it relates global dimensions to the quantum Weyl denominator.
It also contains a definition of the quantum Lie superfactorial, for all cases of Dynkin diagrams.
A few examples and applications are given in the same section, and possible extensions of the present work are mentioned.

Section \ref{McKayCorrespondence} is devoted to an operational description of the correspondence between periodic quivers for Lie groups and the theory of module-categories associated with fusion categories ${\mathcal A}_k(SU(2))$.

Section \ref{sec:periodicquiver} gives a periodic quiver representation of the Weyl vector (or Weyl element), using the correspondence described in section \ref{McKayCorrespondence}.
The analysis is extended to the non simply-laced cases.
A  detailed description is given for the cases $E_6$ and $F_4$, mostly to illustrate the general considerations of the previous section. We obtain, using this approach,  the collection of scalar products needed for the determination of the Weyl denominator. The same technique is  used for the other exceptional cases.  This is also done for several examples taken from the classical series $A_r$, $B_r$, $C_r$, $D_r$.  
In all these cases, the classical Weyl denominator is  obtained as a product of factorials over the exponents, when the Dynkin diagram is simply-laced, with a trivial modification when it is not.
A direct evaluation of this quantity is carried out in the same section for $A_r$, with $r$ generic, using the known $SU(2)$ fusion coefficients,  but this is done for the series $B_r$, $C_r$, $D_r$ (and for $A_r$ again) in the Appendix, without using the quantum McKay correspondence.

\bigskip

{\sl Remark}: the global dimension of $G$ at level $k$ is formally equal to  the square inverse of the Chern-Simon partition function \cite{WittenCS},  for connections on $G$-principal bundles on the sphere $S^3$.  The topics discussed in this paper are therefore also related to differential geometry and to the theory of quantum invariants for knots.

\section{Global dimensions and the Weyl denominator}
\label{sec: globaldim}
\subsection{Main results}
\label{sec: mainresult}
If $n$ is a positive integer and $q$ is complex, or more generally if $n$ an arbitrary real number and $q$ a formal variable\footnote{For $q$ complex, the usual convention is to set $q^n= exp(n \, Log(q))$ while choosing the standard principal value of  $Log$, \ie a branch cut along $(- \infty, 0]$.
While studying non simply-laced cases we shall meet situations where $n$ is a positive fraction with denominator $2$ or $3$; one could define ``v-numbers'' and set $q=v^2$ or $v^3$, but
for the applications that we have in mind, $q=exp(\frac{i \pi}{k+N^\vee})$, and no confusion will arise. }, 
there are two standard conventions for the definition of q-numbers (``quantum numbers''). 
The one\footnote{\label{foot: qnumberconvention}The other standard convention, used for instance in \cite{Mathematica}, defines $q$-numbers as $\frac{1 - q^{n}}{1-q}$.
There are also  two standard conventions for quantum factorials, the one used in the previous reference differs from ours by a prefactor.}  
that we use in this paper is the following:
$$ [n]_q = \frac{q^n - q^{-n}}{q-q^{-1}}$$
We sometimes omit the surrounding square bracket, setting $n_q = ~ [n]_q$.

\bigskip
\noindent
We shall need the following definitions:
\label{sec: def1}
{\definition{}Let $G$ be a simply-laced complex simple Lie group and let $q$ be a complex number. Call ${\mathcal Exp}$ the multiset of exponents of $G$. We define the 
quantum Lie superfactorial of type $G$:
$$ \mathrm{\sfac}_{G}[q]= \prod_{s \in {\mathcal Exp}} \,  [s]!_q $$
}

\noindent
Here $[s]!_q$ denotes\footnote{see footnote \ref{foot: qnumberconvention}.} the quantum factorial of $s$, \ie 
$[s] !_q = \prod_{n=1}^{n=s} [n]_q$.

\bigskip
\noindent
For non-simply-laced cases ($G$ not of type ADE), our definition of the quantum Lie superfactorial functions involves appropriate multiplicative factors that are given 
at the end of section \ref{superfactorial}.

\bigskip
\noindent
For fusion categories of the kind studied in this paper,  the global dimension is given by:
{\theorem{} Let $G$ a simple Lie group with rank $r$, Coxeter number $N$, dual Coxeter number $N^\vee$, and let $k$ be a positive integer.  The global dimension of the fusion category ${\mathcal A}_k(G)$ is given  by 
$$ \vert {\mathcal A}_k(G)\vert \doteq \sum_n \, \mu(n)^2 = 
\frac{(N^\vee+k)^r}{2^{r N} \, \Delta \,  {\left( \sin  \dfrac{\pi}{N^\vee+k}\right)}^{r N} \, \left( \mathrm{\sfac}_G [q]\right)^2}
$$

where $\mathrm{\sfac}_G [q]$ is the quantum Lie superfactorial of type $G$, $\Delta$ is the determinant of the fundamental quadratic form\footnote{\label{foot: symmetrizingcoefficients}Calling $L$ the Cartan matrix, $D= \mathrm{diag}(\delta_s)$ the diagonal symmetrizing matrix normalized in such a way that its largest entries are equal to $1$ (the long roots  have norm square $2$),  and  $K = L^{-1} D$, then
$\Delta = det(K)$ is the inverse of the index of the sublattice of long roots in the weight lattice.
For simply-laced cases, $\Delta = det(L)^{-1}$.}, and $q = \exp{i \pi/(N^\vee+k)}$.
The symbol $\doteq$ stands for ``is defined as'' and the sum of the square of the quantum dimensions $\mu(n)$ runs over the simple objects denoted $n$ (they are finitely many).}
\bigskip

To say that we have a module-category ${\mathcal E}$ associated with a monoidal category ${\mathcal A}$ amounts to say  \cite{Ostrik} that we are given a  
monoidal functor from ${\mathcal A}$ to the category of endofunctors of an abelian category ${\mathcal E}$. 
We then suppose that the fusion category ${\mathcal A}_k(G)$ is given,  and we consider module-categories ${\mathcal E}_k(G)$ associated\footnote{It will be convenient to say that ${\mathcal E}={\mathcal E}_k(G)$ is a module-category of type $G$.} with it.
As before, our purpose is to calculate their global dimensions.
The simplest situation occurs when  the module-category ${\mathcal E}_k(G)$ is associated with the conformal embedding of $Lie(G)$, at level $k$, into another Lie algebra $Lie(J)$, at level $1$. 
The definition of conformal embeddings belongs to the lore of affine Lie algebras but we re-write it here in a way that uses only the properties of finite dimensional Lie algebras.

{\definition{} Let $G$ and $J$ be simple Lie groups, and $k, \ell$ be positive integers.  There is a conformal embedding of $G$ at level $k$, in $J$ at level $\ell$, if the following three conditions are satisfied: 1) There is an embedding of  the Lie algebra $Lie(G)$ into $Lie(J)$, 2) The Dynkin index of the embedding is equal to $k/\ell$, 3) 
If $N^\vee_G$ and $N^\vee_J$ are the dual Coxeter numbers of $G$ and $J$, the following equality  holds: $\frac{dim(G) \times  k }{k + N^\vee_G} = \frac{dim(J) \times \ell}{\ell + N^\vee_J}.$
We call  $c$ the common value of the last two expressions. When\footnote{This is what will be assumed from now on: we set $\ell=1$.} the integer $\ell$ is not specified, it is understood that $\ell=1$.
If $G$ is semi-simple, not simple,  the same definition holds,  but now $k=(k_i)$ is a multi-index, the quantity $c_i$ is defined for each simple component $G_i$ of $G$ and the equality of central charges (condition 3) should hold for $c=\sum c_i$.
}

\smallskip
Conformal embeddings have been classified, see \cite{BaisBouwknegt, KacWakimoto, SchellekensWarner}, and for a given Lie group $G$ at level $k$, they all happen when $\ell = 1$. This justifies the fact of taking $\ell =1$ in the sequel.

\smallskip
Starting with such a conformal embedding $(Lie(G),k) \subset (Lie(J),\ell=1)$, a known construction, recalled in \cite{CordobaRC},  allows one to determine a particular module-category ${\mathcal E}$ that, in some sense, ``measures'' the embedding.

 \label{formulaforE}
 
 {\theorem {} Consider a conformal embedding of the simple or semi-simple Lie group $G$, at level $k$, in the simple Lie group $J$, at level $1$.  This embedding is associated with a module-category 
  ${\mathcal E}={\mathcal E}_k(G)$, with an action of the fusion category ${\mathcal A}={\mathcal A}_k(G)$. The global dimension of ${\mathcal E}$  is given by:
  $$  \vert {\mathcal E} \vert = \sqrt {\vert {\mathcal A}\vert \times   \vert {\mathcal J} \vert } $$
  where   $\vert {\mathcal J} \vert $ is the global dimension of the fusion category  ${\mathcal J}= {\mathcal A}_1(J)$.
} 

 \subsubsection*{About module-categories associated with conformal embeddings}
 We refer to \cite{CordobaRC} for a detailed discussion and proof of theorem 2. 
 Let us just notice that,  since the right hand side of this equality refers to global dimensions of the fusion categories ${\mathcal A}={\mathcal A}_k(G)$ and  ${\mathcal J}={\mathcal A}_1(J)$, 
it can be evaluated immediately thanks to the general formula (theorem 1) without having to compute the quantum dimensions of the simple objects of ${\mathcal E}$. The result can then be expressed in terms of Lie quantum superfactorial functions.

\bigskip

 Simple objects of ${\mathcal A}_k(G)$ correspond to integrable weights of an affine Lie algebra at level $k$ (the elements of the Weyl alcove) associated with a classical Lie algebra $Lie(G)$. In turn, such integrable weights are completely specified by dominant weights $n$ of the underlying classical Lie algebra $Lie(G)$, restricted by the condition $\langle n,  \theta  \rangle \leq k$, where  $\theta$ is the highest root of $Lie (G)$ and $\langle .. , .. \rangle $ is the fundamental quadratic form.  By a  slight abuse of notation, $n$ will denote simultaneously the weight of the affine Lie algebra, the corresponding weight of the classical Lie algebra, and the corresponding irreducible representation.  We shall call {\sl level of a weight} $n$ the integer $\langle n,  \theta  \rangle$ and say that $n$ is integrable at level $k$ or that ``it exists at level $k$'' if and only if its level is smaller than $k$ or equal to $k$.
 If a representation $n$ exists at level $k$, its quantum dimension $\mu(n)$ at that level is given by the quantum version of the Weyl dimension formula recalled below, with $q=exp(i \pi/(N^\vee+k))$. The classical dimension is recovered by taking $q=1$. This expression is usually obtained from the theory of quantum groups or from the theory of affine Lie algebras, but we take it here as a definition.
 Supposing  $n$  irreducible, we  use the same notation\footnote{On purpose, we use the same notation for integers (say $n$), so that, for $G=SU(2)$, an irreducible representation~$n$ of highest weight $n$ has (classical) dimension $n+1$ and quantum dimension $\mu(n)=[n+1]_q$ at level $k$.} for the representation and  its highest weight.
\[\mu (n)= \prod_{\alpha > 0} \frac{\langle n + \varrho, \alpha\rangle_q}{\langle\varrho, \alpha\rangle_q}\] Here $\varrho$ is the Weyl vector, $\alpha$ runs over the set of positive roots, and  $\langle .. , ..  \rangle_q$ denotes the q-number  $[\langle ..  , .. \rangle]_q$ associated with the inner product between the chosen weights (we do not write the surrounding square bracket).
Using the $n_q=[n]_q$ convention to define the $q$-numbers $\langle\varrho, \alpha\rangle_q$, it is convenient to call\footnote{${\mathcal D}_q$ defined below would read differently if we were using  another convention (see footnote \ref{foot: qnumberconvention}) for $q$-numbers.}:
~{\definition{} The quantum Weyl denominator of a Lie group $G$ at level $k$ is  ${\mathcal D}_q=\prod_{\alpha > 0} \,   \langle\varrho, \alpha\rangle_q$}

\smallskip

The last step of the proof of  {\sl Theorem 1}  relies on the following result whose proof is deferred to section \ref{sec:periodicquiver} and to the appendix.

{\theorem{} The quantum Weyl denominator of $G$ is equal to the quantum superfactorial of type $G$, \ie  ${\mathcal D}_q = \mathrm{\sfac}_G[q]$. The equality holds in the classical case, \ie when $q=1$:
the classical Weyl denominator of $G$ is obtained as ${\mathcal D}=  \mathrm{\sfac}_G[1]$.  If $G$ is simply-laced, ${\mathcal D}$ is equal to the product of factorials of the exponents of $G$. If $G$ is not simply-laced this product is modified by a multiplicative factor $\varpi$ equal to $1/2^r$ for $B_r$, $1/2^{r(r-1)}$ for $C_r$, $1/2^{12}$ for $F_4$, $1/3^3$ for $G_2$.}

\smallskip
Using the quantum McKay correspondence for $SU(2)$, we shall relate the matrices $F_n$ describing the module-action of the fusion ring of the category ${\mathcal A}_k(SU(2))$ on the Grothendieck group of a module-category ${\mathcal G}$ characterized by a simply-laced Dynkin diagram $G$ (the index $n$ of $F_n$ refers to the simple objects of ${\mathcal A}_k(SU(2))$ and the entries $(F_n)_{ab}$ of these matrices are indexed by the simple objects of ${\mathcal G}$) to the classical Weyl denominator of type $G$, and therefore also to the corresponding superfactorials.  In particular we have the following ``fusion product formula'', for ${\mathcal D}$: 

{\theorem{} Calling $F_n$ be the fusion matrices describing the module-action of the fusion ring of the category ${\mathcal A}_k(SU(2))$ on the Grothendieck group of a module-category ${\mathcal G}$ characterized by a simply-laced Dynkin diagram $G$, and ${\mathcal D}$ the classical Weyl denominator of the (complex, simply-connected) associated Lie group $G$, of rank $r$ and Coxeter number $N$, we have the following identity, with $\varpi = 1$: 
\[
{\mathcal D}^2 \;= \varpi^2  \times \;\prod_{n=1}^{N} \prod_{b=1}^r \sum_{a=1}^r  \; (F_n+F_{n-1})_{a,b}
\]
With an appropriate redefinition of the matrices $F_n$ (see the text for details), the same equality, with $\varpi$ as in the previous theorem, is also valid when $G$ is not simply-laced. 
}

\subsection{From global dimensions to quantum superfactorials}
\label{sec:globaldimforfusion}

\subsubsection*{Global dimensions and the quantum Weyl denominator}
\label{sec: Smatrix}

\smallskip
\noindent
Action of the group $SL(2,\ZZ)$ on the vector space spanned by the simple objects $m,n\ldots$ of the modular category ${\mathcal A}_k(G)$ is described by unitary matrices $S$ and $T$ representing the two generators $\tau \mapsto -1/\tau$ and $\tau \mapsto \tau +1$ of the modular group. 
This representation,  known by Hurwitz  long ago \cite{Hurwitz} for $G=SU(2)$, is given, for arbitrary simple Lie groups,  by the Kac-Peterson formulae \cite{Kac-Peterson}.  
  A  simple manipulation -- see \cite{CordobaRC} -- of the formula giving $S$ allows one to relate as follows its matrix element $S_{1,1}$ to the quantum Weyl denominator $ {\mathcal D}_q$ (although sometimes written differently, see for instance \cite{YellowBook}, this is certainly well known):
$$
S_{1,1} = \frac {2^{\frac{r N}{2}} \, \sqrt{\Delta}}
{{(N^\vee+k)}^{r/2}} \;  {\left( \sin \frac{\pi}{N^\vee+k} \right) }^{\frac{r N}{2}} \; {\mathcal D}_q$$
Straightforward  manipulations on the quantum Weyl formula lead to the following equality, which is sometimes used, in other contexts, to define quantum dimensions themselves:
$$ \mu(n) = S_{n,1}/S_{1,1}$$
Using unitarity of the $S$ matrix, the global dimension $\vert {\mathcal A}_k(G) \vert$, defined as sum of squares of the q-dimensions of its simple objects, can be written as:
$$\vert {\mathcal A}_k(G)\vert \doteq \sum_n \, \mu(n)^2  = \frac{1}{S_{1,1}^2}$$
From the expression for $S_{1,1}$ given above, one gets:
$$ \vert {\mathcal A}_k(G)\vert = \frac{(N^\vee+k)^r}{2^{r N} \, \Delta \,  {\left( \sin  \dfrac{\pi}{N^\vee+k}\right)}^{r N} \, {\mathcal D}_q^2}$$

\subsubsection*{Expression of the (quantum) Lie superfactorials functions $\mathrm{\sfac}_G[q]$}
\label{superfactorial}

For simply-laced $G$ the functions $\mathrm{\sfac}_G[q]$ have been defined in section \ref{sec: def1}.
For the convenience of the reader, the corresponding explicit expressions are gathered below.

\smallskip
\noindent
Simply laced cases (ADE): 
\begin{itemize}
\item
$G=A_{r}\sim SU(r+1)$, ${\mathcal E} = {1,2,\ldots,r}$, and $ \mathrm {\sfac}_q(r) \doteq \mathrm{\sfac}_{A_r}[q]= \prod_{s=1}^{s=r} \,  [s]!_q $.
\item   $G=D_{r}\sim Spin(2r)$,  ${\mathcal E} = {1,3,5,\ldots N-3, N-1; N/2}$, and\footnote{In the $D_r$ case,  $N=2r-2$, and when $r$ is even, $N/2$ appears twice.} $ \mathrm{\sfac}_{D_r}[q]= [N/2]!_q \, \prod_{s, odd=1}^{s=N-1} \,  [s]!_q $.
\item  $G = E_6$ then $ \mathrm{\sfac}_{E_6}[q]=[1]!_q  \, [4]!_q  \, [5]!_q \, [ 7]!_q  \,  [8]!_q  \,  [11]!_q $  
\item  $G = E_7$  then $ \mathrm{\sfac}_{E_7}[q]=[1]!_q  \, [5]!_q   \,  [7]!_q  \,  [9]!_q  \,  [11]!_q \, [13]!_q \, [17]!_q $
\item  $G = E_8$  then   $\mathrm{\sfac}_{E_8}[q]=[1]!_q  \, [7]!_q \,  [11]!_q \,  [13]!_q \, [17]!_q \,  [19]!_q \,  [23]!_q \,  [29]!_q$
  \end{itemize}

\smallskip  
For non simply-laced cases, the actual calculation of ${\mathcal D}_q$ (see section \ref{sec:periodicquiver}) justifies the definitions of $\mathrm{\sfac}_G[q]$ given below.
Notice that one indeed needs to modify a naive definition by the introduction, into the $q$-numbers entering the $q$-factorial of exponents, of the scaling coefficients  $\delta_s$  defined in footnote \ref{foot: symmetrizingcoefficients}.

\smallskip
\noindent
Non simply-laced cases: 
\begin{itemize}
\item $G=B_{r} \sim Spin(2r+1)$, then ${\mathcal E} = {1,3,5,\ldots ,2r-1}$ and $\mathrm{\sfac}_{B_r}[q]= \prod_{s \in {\mathcal E}} \,  \widetilde{[s]!_q}$ where \newline
$ \widetilde{[s]!_q}  = [\frac{s}{2}]_q \,  [s-1]_q  \,  [s-2]_q \ldots [3]_q  \, [2]_q \,  [1]_q,  \quad \text{and} \quad   \widetilde{[1]!_q} =[1/2]_q$ 
\item $G=C_{r} \sim Sp(2r)$, then  ${\mathcal E} = {1,3,5,\ldots ,2r-1}$ and $ \mathrm{\sfac}_{C_r}[q]= \prod_{s \in {\mathcal E}} \, \widetilde{[s]!_q}$ where \newline
 $\widetilde{[s]!_q}  = \left[\frac{s}{2}\right]_q \, \left[\frac{s-1}{2}\right]_q \ldots
 \, \left[\frac{s-\frac{s-3}{2}}{2}\right]_q \,  \left[{s-\frac{s-1}{2}}\right]_q \, \left[\frac{s-\frac{s+1}{2}}{2}\right]_q  \ldots \left[\frac{2}{2}\right]_q \,   \left[\frac{1}{2}\right]_q$
 $ \text{and} \quad   \widetilde{\left[1\right]!_q} =1$.
\item $G=F_4$, then  ${\mathcal E} = {1,5, 7, 11}$ and\\
  $ \mathrm{\sfac}_{F_4}[q]=\left[\frac{1}{2}\right]_q^2 1_q^3 \left[\frac{3}{2}\right]_q [2]_q^3
   \left[\frac{5}{2}\right]_q^2 [3]_q^3 \left[\frac{7}{2}\right]_q [4]_q^2
   \left[\frac{9}{2}\right]_q [5]_q^2 \left[\frac{11}{2}\right]_q [6]_q [7]_q
   [8]_q $
\item $G=G_2$, then  ${\mathcal E} = {1,5}$ and  $ \mathrm{\sfac}_{G_2}[q]=  \left[\frac{5}{3}\right]_q \left[\frac{4}{3}\right]_q 3_q \,   2_q \,  1_q \,   \left[\frac{1}{3}\right]_q$
\end{itemize}

\bigskip

{\sl Relation with the Barnes G-function}.
\noindent
When the argument is an integer, the quantum superfactorial $ \mathrm {\sfac}_q(r)$, that we associate with $A_r$,  differs from the (shifted) quantum Barnes G-function by a pre-factor. The latter
 disappears if one uses the alternative convention for the definition of $q$-numbers, and the corresponding quantum superfactorials, alluded to in footnote \ref{foot: qnumberconvention}. We refer to \cite{CordobaRC} for a discussion of the relations with the quantum Barnes G-function. 
Note that the (classical, \ie $q \mapsto 1$) Lie superfactorial of type $A_r \sim SU(r+1)$  coincides with the usual superfactorial $\mathrm{\sfac}(r)$ and therefore with the value of the  (classical) Barnes function $\text{G}(r+2)$. More generally, the classical limit of $\mathrm{\sfac}_G$ gives interesting sequences -- or particular integers -- that have been added to \cite{OEIS}.

\subsection{Examples and applications}
\label{sec:examples}

\subsubsection*{Examples for \texorpdfstring{${\mathcal A}_k(G)$}{AkG} and applications}
The main theorem giving $\vert {\mathcal A}_k(G)\vert$ leads, for every choice of $G$, to fully explicit formulae that can be expressed in terms of elementary functions, since 
$[n]_q={\sin \left(\frac{\pi  n}{k+N}\right)}/{\sin \left(\frac{\pi }{k+N}\right)}$.
For convenience, we remind the reader the values of Coxeter numbers $N$, dual Coxeter numbers $N^\vee$, and  ``long indices''  $\Delta^{-1}$ for all Lie groups (the quantity $\Delta$, used in the text, is the determinant of the fundamental quadratic form,  \ie the inverse of the long index).

\smallskip
\begin{center}
$
\begin{array}{cccccccccc}
{} & A_r & B_r& C_r& D_r & E_6 & E_7& E_8 & F_4& G_2 \\
N : & r+1 &2r & 2r &2r-2 & 12& 18 & 30 & 12& 6 \\
N^\vee :  &  r+1&2 r -1 &r+1 &2r-2 &12 &18 &30 & 9& 4 \\
\Delta^{-1}: & r+1 & 4 & 2^r & 4 & 3 & 2 & 1 & 4 & 3
\end{array}
$
\end{center}

Using these values together with the definition of the quantum Lie  superfactorial functions, one finds, for instance, that the global dimension of $A_{N-1}\sim SU(N)$ at level $k$ reads:
$$
\vert {\mathcal A}_k (SU(N)) \vert =\frac{N (k+N)^{N-1}}{2^{N (N-1)} \prod _{s=1}^{N-1} \sin ^{2
   (N-s)}\left(\frac{\pi  s}{k+N}\right)}
 $$
whereas, in the case of $G_2$, for instance, it is:
$$
\vert {\mathcal A}_k (G_2) \vert=
\frac{3 (k+4)^2}{2^{12} \sin ^2\left(\frac{\pi }{3 (k+4)}\right) \sin
   ^2\left(\frac{\pi }{k+4}\right) \sin ^2\left(\frac{4 \pi }{3 (k+4)}\right) \sin
   ^2\left(\frac{5 \pi }{3 (k+4)}\right) \sin ^2\left(\frac{2 \pi }{k+4}\right) \sin
   ^2\left(\frac{3 \pi }{k+4}\right)}
   $$
 These expressions can be used to study various kinds of limits. 
If we remember that $r+r N$  gives the dimension $dim_G$ of the Lie group $G$, we obtain the following ``classical limit'' : 
  $$ \vert {\mathcal A}_k(G)\vert \,  \underset{{{k\rightarrow \infty}} }{\sim}\, \frac{ k^{dim_G}}{(2\pi)^{rN} \, \Delta \,  \left({\mathrm{\sfac}_G}\right)^2}
\;  {\text {, in particular, }} \;  \vert {\mathcal A}_k(SU(N))\vert \,  \underset{{{k\rightarrow \infty}} }{\sim}\, \frac{N\times k^{N^2-1}}{(2\pi)^{N(N-1)}  \,  \, \left({\mathrm{\sfac}(N-1)}\right)^2}
  $$

\noindent
The level-rank duality property for Lie groups of type $A_{N-1} \sim SU(N)$  was observed  in \cite{JimboMiwa}; for the coefficient $S_{1,1}$ of the modular matrix $S$ it implies:
$\sqrt{N} \;  S_{1,1}[A_{N-1},k] = \sqrt{k} \; S_{1,1}[A_{k-1},N] $.
In terms of global dimensions we have therefore
$
k{\vert {\mathcal A}_k (SU(N)) \vert} = N {\vert {\mathcal A}_N (SU(k)) \vert}.
$
This equality can also be used to get simple enough expressions for global dimensions of $SU(N)$ at level $k$, when the rank $N-1$ is large and the level $k$ reasonably small. 
Using the previous asymptotic expression of $\vert {\mathcal A}_k (SU(N)) \vert$ for large $k$, the  duality relation for $SU(N)$ gives immediately another asymptotic expression when the level $k$ is fixed and the rank $r=N-1$ goes to infinity:
$$
\vert {\mathcal A}_k(A_{r}) \vert  \,  \underset{{{r \rightarrow \infty}} }{\sim}\, \frac{1}{(2 \pi )^{k(k-1)} \,   \left[{\mathrm{\sfac}(N-1)}\right]^2} \times r^{(k^2)}
$$

\subsubsection*{Examples for \texorpdfstring{${\mathcal E}$}{E} and applications}

The list of conformal embeddings has been known for more than twenty years, see \cite{BaisBouwknegt, KacWakimoto, SchellekensWarner}.
For the sake of illustration we only remind the reader the list of conformal embedding for $G=SU(N)$ Lie groups.
There are three regular series and a few sporadic cases:
 
 \smallskip
 \noindent 
{\sl Regular series}. The following conformal embeddings are respectively called antisymmetric, adjoint and symmetric (take $G=SU(N)$ in all cases):
 $$ 
 \begin{array}{c|c|c|c}
k &   k=N-2, &  k=N, & k=N+2, \\
{} &       N \geq 4 & N \geq 3 & N\geq 2 \\
J & SU({N(N-1)}/{2}) & Spin(N^2-1) & SU(N(N+1)/2)  \\
 \end{array}
$$ 
 
\noindent
{\sl Sporadic cases}.
$$
 \begin{array}{c|cc|cc|c|c|cc|c}
 G & SU(2)  & {} &  SU(3) & {}  & SU(4) & SU(6) & SU(8) & & SU(9) \\
 \hline
 k & 10 & 28 &9& 21&  8 &  6 & 1 & 10 &  1\\
J & Spin(5) & G_2  & E_6 & E_7  & Spin(20)  & Sp(10) & E_7 & Spin(70)  & E_8 \\
 \end{array}
$$ 

\smallskip
\noindent
For module-categories ${\mathcal E}$ associated with  symmetric or antisymmetric regular conformal embeddings of $SU(N)$ into $J = SU(p)$ with $p=N(N\pm1)/2$, the calculation of $\vert {\mathcal E} \vert$ using theorem 2 is particularly simple, since, for any integer $p>1$, $\vert{\mathcal A}_1(SU(p))\vert = p$. Indeed, besides the trivial representation,  the $p-1$ fundamental irreps of $SU(p)$ all exist at level $1$ and they have quantum dimension~$1$. Notice that for $G=SU(2)$ the only regular embedding occurs at level $4$, it is symmetric and is associated with a module described by the fusion graph $D_4$ (the usual Dynkin diagram of $Spin(8)$). Still with $G=SU(2)$ we have also two sporadic embeddings, into $Spin(5)$ and $G_2$, respectively described (and usually denoted) by the fusion diagrams $E_6$ and $E_8$. The global dimension of their associated modules can be determined from theorem 2 using $\vert {\mathcal A}_{1}(SU(3))\vert =3$,   $\vert {\mathcal A}_{1}(Spin(5))\vert =4$ and $\vert {\mathcal A}_{1}(G_2)\vert =\frac{1}{2} (5 + \sqrt{5})$ together with $\vert {\mathcal A}_k(SU(2)) \vert = ({1}/{2}) {(k+2)}/{\sin^2({\pi}/{(k+2)}})$. In those particularly easy cases, the results could be directly obtained by summing  the squares of the quantum dimensions of the simple objects. Explicit results for other examples can be found in \cite{CordobaRC}.
Before ending this section, let us mention that the adjoint conformal embedding exists for all Lie groups, it is an embedding  of $G$ at level $k=N^\vee$, its dual Coxeter number, into $Spin(dim(G))$.

\subsubsection*{Higher sums}

The global dimension that was defined by summing squares of 
quantum dimensions $\mu(n)$ over the simple objects $n$ of a category $\mathcal{A}_k(G)$
appears as a particular case, for $s=-2$,  of a generalized Riemann Zeta function defined as $\zeta_G(s,k) = \sum_n \mu(n)^{-s}$,
and simple formulae should certainly exist.
For instance, taking  $G=SU(2)$, we know that $\vert \mathcal{A}_k(SU(2)) \vert = \zeta_{SU(2)} (-2,k)=\frac{1}{2} (k+2) \csc ^2\left(\frac{\pi }{k+2}\right)$, but we found some experimental evidence that 
 for $p$ integer,  both even and positive, the quantity $\zeta_{SU(2)} (-p,k)$ is given, for levels $k > p/2 -1$, by the expression
$(k+2) \,  4^{-p/2} \binom{p}{\frac{p}{2}} \sin ^{-p}\left(\frac{\pi}{k+2}\right)$.
When $p$ is even but negative, more involved expressions can also be found. 
We hope that the perspective of proving or generalizing these observations may trigger the interest of the reader.

\section{Periodic quivers and quantum McKay correspondence}
\label{McKayCorrespondence}

The ADE correspondence between module-categories of type $SU(2)$ and simply-laced Dynkin diagrams was  first obtained by theoretical physicists in the framework of conformal field theories (CFT) (classification of modular invariant partition functions for the WZW models of type $SU(2)$, \cite{CIZ}). It was set in a categorical framework by \cite{Ostrik, KirilovOstrik}. 
In plain terms, the diagrams encoding the action of the fundamental irreducible representation of $SU(2)$ at level $k$ (which is classically $2$-dimensional) on the various module-categories existing at that level,  are the Dynkin diagrams describing the simply-laced simple Lie groups with Coxeter number $k+2$. This is nowadays a well-known result.
What is maybe not so well known is that, at a deeper level, there is a correspondence between fusion coefficients of the $SU(2)$ module-category described by a Dynkin diagram $G$
and the numbers obtained by calculating the scalar products between  fundamental weights and (all) the roots of the associated simply-laced Lie group.
Although studied in another context, this  was described in \cite{Dorey:CoxeterElement}, and, while expressed in a different manner, it was  independently stated in \cite{Ocneanu:Bariloche} where one can also find a discussion of what happens when $SU(2)$ itself is replaced by an arbitrary simple Lie group -- a generalization that we do not need in the present article.
In the non-ADE cases, and although one cannot associate $SU(2)$ module-categories to non simply-laced Dynkin diagrams, one can nevertheless use the 
action of the fusion ring on the modules defined by the chosen Dynkin diagrams, and, with a simple modification of the rules,
 still obtain a correspondence between the coefficients describing the module structures and the families of scalar products already mentioned (one has only to introduce scaling coefficients in appropriate places). 
As the property relating fusion coefficients and scalar products is not well documented, and since we shall use it in section \ref{sec:periodicquiver}, we shall 
review it in more detail and sketch in section \ref{sec:fusiontoscalarproducts} the origin and proof of this relation.

\subsection{The bipartite Coxeter element and the set ${\mathcal R}$}
Being bipartite we can color in black or white the nodes of a chosen Dynkin diagram  and consider separately the product of black or white simple Weyl reflections associated with the corresponding nodes. The black nodes  define a commuting set of involutions, so their product is also an involution; we have the same thing for the white nodes.
The product of these two involutions is called a bipartite Coxeter element (call it $c$). One then considers the action of this Coxeter element on the set of roots and decompose the latter into orbits.  Still calling $r$ the rank and $N$ the Coxeter number, one shows that these $r$ orbits have $N$ elements each (use the fact that $c^N=1$). This observation may be traced back to \cite{Bourbaki:groups, Kostant:SU2, Steinberg} and was used more recently, see  \cite{KirillovThind} to attach to every simple Lie group a particular periodic quiver whose vertices are labelled by roots. 
This combinatorial structure was already present in reference \cite{Ocneanu:Bariloche}. 
 In order to better appreciate the relation between this construction and the theory of module-categories for $SU(2)$,  one may proceed as follows.
For a simply-laced Lie group with Dynkin diagram $G$ and adjacency matrix also called $G$,  Coxeter number $N$ and rank $r$, let us call $\mathcal R$ the set $(\ZZ \times_{\ZZ_2} G)/\ZZ_{2N}$. Its $N \times r$ elements can be displayed as the pattern obtained by considering only the dark tiles of a periodic rectangular checkerboard of width $r$, height $2N$, periodic in the vertical direction. Using coordinates, $\mathcal R$ is the set of ordered pairs $(n,b)$, where $n \in {1,2,\ldots, 2 N}$ and $b \in {1,2,\ldots r}$ such that $n+b$ is even.
Elements of the vector space $\CC^{Nr}$ of complex valued functions on the set  $\mathcal R$ can be displayed as rectangular matrices with $2 N$ lines and $r$ columns, where we only consider those entries that are such that $n+b$ is even. For short, we  say that such entries are the ``even positions'' of the underlying matrix.  We shall  display these functions as rectangular arrays of numbers located at even positions and we shall not draw the surrounding brackets to keep in mind the distinction between such tables and usual matrices. 
As the cardinality of the set of roots is also $N \times r$, it is obvious that $\mathcal R$ is in bijection, in many ways, with the set of roots of the chosen Lie group. 
We shall see how $SU(2)$ fusion matrices describing an appropriate module-category establish between these two sets a special one-to-one correspondence enjoying nice properties. Using this particular bijection one considers the location of a root as a particular point of $\mathcal R$ and associates  the columns of $\mathcal R$ with the orbits of the bipartite Coxeter element $c$.  In other words roots can be obtained as particular elements of the vector space $\CC^{Nr}$, namely Dirac measures centered on the points of $\mathcal R$.

\subsection{Periodic essential matrices (simply-laced cases)}
\label{sec:simplylaced}

The fusion category ${\mathcal A}_k(SU(2))$ possesses  $k+1$ simple objects ${(\und n)} \in {(\0),(\1),(\2),\ldots, (\und k)}$. 
The Grothendieck ring (the fusion ring)  has one algebraic generator $\sigma = (\1)$ and the unit is $(\0)$. The other linear generators are obtained by the  Chebyshev  recurrence  ${(\und n)=(\und{n-1}) \, \sigma - (\und{n-2})}$.  Moreover, we have one relation: 
$(\und{k+1})=0$ which implies $(\und{k-1})=(\und{k})\,\sigma$. These well known results rely on the fact that the same recurrence  (but not the extra relation $(\und{k+1})=0$) holds for the representation theory of $SU(2)$ itself. 

In the previous sections, the symbol ${\mathcal E}$ was used to denote some module-category of type $G$ (a Lie group), \ie associated with ${\mathcal A}_k(G)$, but in the present section we shall, on purpose, always denote the module-categories of type $SU(2)$ by the symbol ${\mathcal G}$.
In other words, ${\mathcal G} = {\mathcal E}_{N-2}(SU(2))$. 
The chosen Dynkin diagram $G$ is called the fusion diagram of ${\mathcal G}$ because its adjacency matrix (also called $G$) describes the action of the generator $\sigma$ of the fusion ring of $SU(2)$ on the simple objects of ${\mathcal G}$: $\sigma \times a = \sum_b \, G_{a,b} \, b$.
Indecomposable module-categories ${\mathcal G}$ of type $SU(2)$, \ie associated with ${\mathcal A}_k(SU(2))$, are classified by simply-laced Dynkin diagrams. 
Because of this, to every vertex $a$ of a simply-laced\footnote{dual Coxeter number  and Coxeter number coincide: $N^\vee=N$}  Dynkin diagram $G$ with $r$ vertices and Coxeter number $N$ one can associate {\sl both} a fundamental weight $\omega_a$ of a Lie group $G$, and a simple object $a$ of an  $SU(2)$ module-category ${\mathcal G}$ at level $k = N -2$ (this $2$ stands for the dual Coxeter number of $A_1 \sim SU(2)$).  When ${\mathcal G} = {\mathcal A}_k(SU(2))$, the Dynkin diagram is $A_{k+1}$, and the choice of a vertex  determines either a fundamental weight of $SU(k+2)$ or  a simple object of ${\mathcal A}_k(SU(2))$ belonging to the ordered set  $(\0),(\1),(\2),\ldots, (\und k)$. This vertex can also be referred to by its position (notice the shift by $1$) in the list $(1,2,\ldots, k+1)$.

The adjacency matrix of the chosen Dynkin diagram is symmetric since we assumed that it is simply-laced (the Cartan matrix is  ${2}  - G$). 
We then consider the $r \times r$ matrices $F_n=F_{(\und{n-1})}$ defined by the following Chebyshev recurrence relation, that implies $F_{2}=G$.
 $$F_{n}= F_{n-1} . G - F_{n-2}, \qquad F_{0} = 0,  \qquad F_{1} = \one$$  
 These matrices, that we call\footnote{Some authors reserve this name to the case where ${\mathcal G}={\mathcal A}_k(SU(2))$, \ie $G=A_{k+1}$.} ``fusion matrices'', for short, are sometimes called ``nimreps'' or ``annular matrices'' in CFT. They encode the module-action of the fusion ring of $SU(2)$ at level $k$ on the Grothendieck  group of ${\mathcal G}$, which is automatically a $\ZZ_+$ module (its structure constants $(F_n)_{ab}$, labelled by simple objects $a,b$ of ${\mathcal G}$, and defined by $(\und n) a = \sum_b  (F_{(\und n)})_{ab}\; b$ are non negative integers).  The matrices themselves are labelled by the simple objects of ${\mathcal A}_k(SU(2))$, \ie by integers belonging to the prescribed range, but if we take instead $n \in \mathbb{N}$ we notice that the sequence $F_n$ is periodic with period $2 N$, where $N=k+2$ is the Coxeter number. 
The first $k+1$ members of the sequence have non-negative integer matrix elements, but  $F_0=F_N=0$, and the next $k+1$ members have non-positive matrix elements since they obey $F_{N+m} = - F_{N-m}$ for all positive integers $m$; then, by periodicity  $F_{2N}= 0$, \etc
From the family of square matrices $F_n$ one defines the rectangular matrices $\mathfrak{r}_a = (\mathfrak{r}_a)_{nb} = (F_n)_{ab}$, where $a,b \in G$. 
By taking $1\leq n \leq N-1$ we obtain $r$  rectangular matrices $\mathfrak{r}_a$ (one for each vertex of the Dynkin diagram $G$) with $N-1$ lines and $r$ columns; they are called intertwiners  in conformal field theory, or ``essential matrices'', and their entries are non-negative integers.
 If, instead,  we allow the index $n$ to run from  $0$ to $2N-1$, we obtain $r$ rectangular matrices of size $r \times 2N$  that have matrix elements of both signs, but these matrices should be considered as periodic in the vertical direction because of the the periodicity of matrices $F_n$;  moreover the first line is null. In a similar way, if we take $1 \leq n \leq 2N$, the line indexed by $N$ is null.

\smallskip
ADE Dynkin diagrams being bipartite we declare the first vertex $v$ of a long branch of the chosen diagram $G$ to be even ($\partial \, v =0$), and grade the other vertices accordingly.  This defines a $\ZZ_2$ grading $a \in G \rightarrow \partial{a} \in \ZZ_2$ on the vertices. 
  In particular, for $G=A_{N-1}$, we set $\partial{(\0)} = \partial{1}=0$,  $\partial{(\1)}=\partial{2}=1$, and more generally $\partial{(\und s)} = s  \, {mod} \; 2$, 
  so that the grading is also compatible with the module action.
  We now define, for any vertex $a$ of $G$ the rectangular matrices $\mathfrak{e}_a$  of size $r \times 2N$ and 
   matrix elements $(\mathfrak{e}_a)_{nb}$ with $b \in G$ and $n \in {1,2,\ldots, 2 N}$:
 \begin{eqnarray*}
\text {If} \quad  \partial {a} = 0,  \quad   (\mathfrak{e}_a)_{nb} &=& (F_n)_{ab} 
\\
\text {If} \quad  \partial {a} = 1,  \quad   (\mathfrak{e}_a)_{nb} &=& (F_{n-1})_{ab} 
\end{eqnarray*}
In the same way, we define  rectangular matrices $\mathfrak{e}^{+}_a$ of size $r \times N$ by just restricting the range of $n$ indices  to $n \in {1,2,\ldots, N}$.
Notice that the first line (resp. the last line) of $\mathfrak{e}^{+}_a$ is null when $a$ is odd (resp. even).

For an arbitrary weight $\omega = \sum_a \lambda_a \, \omega_a$, of the Lie group $G$,  with $\lambda_a \in \ZZ$ if $\omega$ is integral, we now set
 $\mathfrak{e}(\omega) = \sum_a \lambda_a \,\mathfrak{e}(\omega_a)$ and  $\mathfrak{e}^{+}(\omega) = \sum_a \lambda_a \,\mathfrak{e}^{+}(\omega_a)$, with
 $\mathfrak{e}(\omega_a)= \mathfrak{e}_a$ and  $\mathfrak{e}^{+}(\omega_a)= \mathfrak{e}^{+}_a$.
  In particular, for the Weyl vector  $\varrho = \sum_a \, \omega_a$, we obtain $\mathfrak{e}^{+}(\varrho) =  \sum_a \, \mathfrak{e}^{+}_a$.
 Because of the horizontal $\ZZ_2$-grading (bipartition of the vertices of $G$) and of the vertical $\ZZ_2$-grading (Chebyshev recurrence formula for $F_n$)  the matrix elements  of $\mathfrak{e}_a$ vanish in position $(n,b)$ when $n$ and $b$ are not both even, or both odd, for all possible choices of $a \in G$.
 Obviously, the table obtained by keeping only those entries where $n+b$ is even can be identified with a function on the set ${\mathcal R}$ defined in the previous subsection.
  So, any weight $\omega$  determines a matrix $\mathfrak{e}(\omega)$, and therefore also a table of the above type, that we still call $\mathfrak{e}(\omega)$, and that we think of as a function on $\mathcal R$.

\subsection{Fusion numbers and scalar products}
\label{sec:fusiontoscalarproducts}
\label{sec:rootstoribbon}

The correspondence between fusion coefficients (of a module-category ${\mathcal G}$ over ${\mathcal A}_k(SU(2))$ described by a simply-laced Dynkin diagram $G$) and scalar products (relative to the weight and root lattices of the Lie group $G$) relies essentially on the observation \cite{Dorey:CoxeterElement, Ocneanu:Bariloche} that both can be explicitly calculated with the same formulae:

The eigenvalues of the (symmetric) matrix $G=F_{(\1)}=F_2$ are $\Delta_1 = 2 \cos( \epsilon_j \, \pi/N)$, where  $\epsilon_j$, ${j \in 1,\ldots, r}$ belong to the multiset of exponents, denoted  ${\mathcal E}xp$. Their set is included into the set of exponents of ${\mathcal A}_k(SU(2))$, described by the Dynkin diagram $A_{k+1}$,  \ie into  ${ 1,\ldots, N-1}$ with $k=N-2$.
From this point of view the exponents of ${\mathcal G}$, or of $G$, can be considered as a particular family of  simple objects of ${\mathcal A}_k(SU(2))$, possibly coming  with  multiplicities.
A more traditional approach to exponents defines them from the eigenvalues $\exp(2 i \pi \epsilon_j /N)$ of some (arbitrary) Coxeter element, but the point of view adopted here, using  the adjacency matrix of a fusion graph (here a Dynkin diagram) and Perron-Frobenius theory, can be generalized to higher situations where $SU(2)$ is itself replaced by a more general Lie group.

Since the $F_n$ matrices can be obtained from $F_2=G$ by a Chebyshev recurrence, they can be expressed as matrix polynomials $F_n=U_{n-1}(G/2)$, where $U_n(z)=\frac{1}{n!} \; \partial^n f[t,z]/ \partial t^n \; \vert_{t=0}$ are the usual U-Chebyshev polynomials with generating function $f[t,z]=\frac{1}{t^2-2 t z+1}$ that are such that $U_n[\cos \theta] = \sin((n+1) \theta) / \sin(\theta)$.
For this reason, the $F_n$'s constitute a commuting family of matrices with eigenvalues $\Delta_{n-1}$ with  $\Delta_n = {\sin((n+1) \epsilon_j  \, \pi / N)}/{\sin( \epsilon_j \pi / N)}$.
Calling $\Psi$ an $r\times r$ unitary diagonalizing matrix for $G$, we write $F_n = \Psi^\dag \, \text{diag}(\Delta_n) \, \Psi$, where $\text{diag}(\Delta_n)$ is a diagonal matrix with the above eigenvalues. 
If $G=A_{k+1}$, the exponents are consecutive integers, and the diagonal matrix $\text{diag}(\Delta_n)$ coincides with the row $n$ of the ``modular matrix'' $S$ of section \ref{sec: Smatrix}, normalized by its first element, \ie $\text{diag}(\Delta_n) = S_{n,1}/S_{1,1}$.  Indeed, in the case of $SU(2)$, the matrix $S$ is given by $S_{m,n}=  \sqrt{{2}/{N }} \sin \left({\pi \, m \, n}/{N }\right)$, and $\Psi$ can be taken as $S$ itself. 
So, if $G=A_{k+1}$,  the previous formula for $F_n$ expresses  the fusion matrices in terms of matrix elements of $S$ alone,  an identity  known as the Verlinde formula \cite{Verlinde}.
$S$ is always unitary and symmetric, and it obeys $S^4=1$, but for $SU(2)$ it is real, and obeys $S^2=1$.
Going back to the general case, we write  $F_n = \Psi^\dag \, \text{diag}(\Delta_n) \, \Psi$, where $\Psi_{\epsilon_j} = (\Psi_{a \epsilon_j})$ is an eigenvector of $G$ for the eigenvalue $\epsilon_j$.
 
For every exponent $\epsilon_j$ of $G$ one may consider a  $2$-dimensional  subspace of $\RR^r$ that is simultaneously left invariant by the two involutions whose product defines the already introduced bipartite Coxeter element $c$. This subspace is also invariant under $c$, where the latter acts as a rotation of angle $2  \pi \epsilon_j /N$. By decomposing the simple roots on these invariant subspaces it is not too difficult to determine explicitly the values of the scalar products between any fundamental weight (say $\omega_a$) and the points (roots) belonging to the orbit 
of  a simple root  (say $\alpha_b$) under  the action of $c$.
The obtained expressions distinguish the cases where $a$ and $b$ refer to even or odd nodes of the Dynkin diagram and the result -- see formula (2.14) of \cite{Dorey:CoxeterElement} -- almost coincides with the one obtained for the fusion coefficients $(F_n)_{ab}$, namely $F_n = \Psi^\dag \, \text{diag}(\Delta_n) \, \Psi$.
In order for the two expressions to coincide exactly, one has to label orbits, not by simple roots, but by another system of orbits representatives for the $c$ action (see formula (2.13) of \cite{Dorey:CoxeterElement}, see also \cite{Kostant:SU2}).  This labeling issue is not important for us, and what matters is the fact that the two lists are the same. 
The expression independently obtained by  \cite{Ocneanu:Bariloche} expresses the scalar products between roots themselves, taken as arbitrary elements of the set ${\mathcal R}$, in terms of fusion coefficients.
Our choice of displaying, rather, the scalar products between fundamental weights and (all) roots is essentially equivalent to the approach\footnote{\label{footnote:preprint} This last reference was mentioned to us after the release of our preprint, it has been added to the revised version of the original manuscript.} of \cite{Dorey:CoxeterElement} but our formulation using the rectangular matrices $\mathfrak{e}_a$, thought of as functions on the table ${\mathcal R}$, may be simpler (remember that our definition of the numbers $(\mathfrak{e}_a)_{nb}$ uses $F_n$ if $\partial a = 0$ but $F_{n-1}$ if  $\partial a = 1$).

This discussion justifies the fact that, once a set of simple roots (or a set of fundamental weights) has been chosen for $G$, one can define the coordinates of a root as a pair $(n,b)$ parametrizing the set ${\mathcal R}$ and obtain the scalar product between some fundamental weight, say $\omega_a$, and the root with ${\mathcal R}$-coordinates $(n,b)$ as the number $(\mathfrak{e}_a)_{nb}$ defined in the previous section in terms of fusion coefficients.  The tables $\mathfrak{e}_a$, immediately obtained from the fusion coefficients $(F_n)_{ab}$, contain enough information to identify all the roots $(n,b)$ and express them -- for instance -- in terms of simple roots, and to study the orbits of $c$. 
Permuting ``even'' and ``odd'' in the definition of  $\mathfrak{e}$ amounts to choose the other orientation for ${\mathcal R}$, \ie to decompose the set of roots into orbits of the Coxeter element $c^{-1}$ rather than of $c$.
We summarize the discussion as follows:

\smallskip  \noindent
\label{theo:theband}  
{\sl Let  $\omega$ be an arbitrary weight associated with the simply-laced Dynkin diagram $G$, of Coxeter number $N=k+2$,  and consider the table  $\mathcal R$, with entries $(\mathfrak{e}(\omega))_{nb}$ defined as previously in terms of fusion coefficients of  the category ${\mathcal A}_k(SU(2))$ acting on a module-category ${\mathcal G}$ also described by the graph $G$. 
Using the partition of the set of roots of the Lie group $G$ into orbits of a bipartite Coxeter element, one can label the $N \times r$ entries of the table $\mathcal R$ with the roots of the Lie group $G$ in such a way that the numbers $(\mathfrak{e}(\omega))_{nb}$ give the scalar products between the weight $\omega$ and the roots of $G$ located at positions $(n,b)$ of the table $\mathcal R$.
}

\smallskip

Notice  that, as roots are weights, and since roots have a norm square equal to $2$, the precise position in $\mathcal R$ of an arbitrary root $\alpha$ can be determined by decomposing this root on the basis of fundamental weights and calculating the matrix $\mathfrak{e}(\alpha)$:  the obtained table will contain only one entry equal to $+2$ and this will give the position of $\alpha$ in $\mathcal R$. Such detailed calculations are not needed for our purpose but for the sake of illustration we shall  display in figure \ref{table:rootsofE6}, section \ref{Example:E6}, the position of the $72$ roots of $E_6$.
We may also notice that since the matrix $F_2=G$ is the adjacency matrix of the chosen Dynkin diagram, the $SU(2)$ Chebyshev recurrence relation for fusion matrices $F_n$ gives us a relation between entries of the tables ${\mathcal R}$,  \ie between scalar products. For simple weights, and therefore for arbitrary weights, it reads immediately: 
$\mathfrak{e}(\omega)_{n,b}  + \mathfrak{e}(\omega)_{n+2,b} =  \sum_{c \; , G_{bc} \neq 0} \mathfrak{e}(\omega)_{n+1,c} $.
In other words, for each empty position $(n+1,b)$ of the table ${\mathcal R}$, \ie $n+1+b$ odd,  the sum of the vertical neighbors must be equal to the sum of the horizontal neighbors (horizontal, in the sense defined by $G$). For instance, with $G=E_6$, and $\omega=\rho$, the Weyl vector, we read 
$\mathfrak{e}(\rho)_{(n=7, b = 3)} + \mathfrak{e}(\rho)_{(n=9, b = 3)} = \mathfrak{e}(\rho)_{(n=8, b=2)} + \mathfrak{e}(\rho)_{(n=8, b=4)}, + \mathfrak{e}(\rho)_{(n= 8, b=6)}$, \ie $11+9= 7+7+6 $.  See the $E_6$ example in
section \ref{sec:periodicquiver}.

\subsection{Periodic essential matrices (non simply-laced cases)}
\label{sec:nonADEessmat}
When the simple Lie group $G$ of rank $r$ is not simply-laced we still call $G$ the non-symmetric  matrix such that $2 - G$ is the Cartan matrix.
$G$ is the adjacency matrix of a directed graph (not a Dynkin diagram) that we associate with the chosen group.
We  still define matrices $r \times r$ matrices $F_n$, with integer matrix elements,  by the same Chebyshev recurrence relation as in the simply-laced case (section \ref{sec:simplylaced}) and still call them fusion matrices. 
The sequence is again periodic, of period $2N$, where $N$ is the Coxeter number of $G$, which now differs from the dual Coxeter number $N^\vee$.  The number of roots, as the number of coroots, is $r \times N$, but roots and coroots now differ.
For Lie algebras $G_2$, $F_4$ and $B_r$ we follow the convention of drawing the short roots to the right, with scaling factors\footnote{See footnote \ref{foot: symmetrizingcoefficients}.}  $\delta_s$  respectively equal to $(1, 1/3)$, 
  $(1, 1, 1/2, 1/2)$, and $(1,1,\ldots, 1, 1/2)$. For the Lie algebra $C_r$ we draw short roots to the left, with scaling factors $(1/2, 1/2, \ldots, 1/2, 1)$. Remember that long roots have norm $2$.
   Using the above drawing conventions, we attribute a $Z_2$ grading $\partial v = 0$ to the leftmost vertex $v$, and grade the other vertices accordingly, in an alternating way. 
 We now can, in exactly the same way as in the simply-laced case,  define tables  $\mathfrak{e}^{+}(\omega)$ or  $\mathfrak{e}(\omega)$ for each weight $\omega$ of $G$.

  For $SU(2)$ at level $k$, the duality $n \rightarrow \overline{n}$, corresponding to complex conjugation of representations, is trivial. 
As it is well known, this rules out the non-ADE Dynkin diagrams as describing possible fusion graphs for module-categories over ${\mathcal A}_k(SU(2))$ since the matrices describing the structure constants of the fusion ring of  the latter, or its action on the associated modules, should be symmetric. 
Nevertheless like in the ADE cases, the matrix elements of the $F_n$'s (or the entries of the $\mathfrak{e}_a$'s) still describe the scalar products between fundamental  weights $\omega_a$ and  the {\sl coroots} of a non simply-laced Lie group associated with the adjacency matrix $G$. In order to describe scalar products between  $\omega_a$ and the {\sl roots} of $G$, one has to introduce independent scaling factors for the different columns of the tables $\mathfrak{e}_a$ -- equivalently, for the orbits of the corresponding bipartite Coxeter element -- equal to the scaling factors relating simple roots and simple coroots. These factors are non-trivial (not equal to $1$) when the corresponding orbit is an orbit of short roots. We shall denote  ${}{\mathfrak{f}^{}_a}$, and more generally ${}{\mathfrak{f}^{}(\omega)}$, for an arbitrary weight $\omega$,  the table --  the function on ${\mathcal R}$ -- encoding this information, \ie the collection of scalar products between $\omega$ and all the roots (for ADE cases, we have $\mathfrak{f}=\mathfrak{e}$). In particular, the Weyl denominator used in sections 1 and 2 of the paper is obtained by multiplying together all the entries of ${}{\mathfrak{f}^{+}(\varrho)}$.

\subsection{From fusion matrices to the Weyl denominator}
\label{sec:weyldenominatorfromfusion}
The scalar product between the Weyl vector $\rho$ of a Lie group associated with the Dynkin diagram  with Cartan matrix ${2}\one  - G$ and the coroot $\beta$  localized at the position $(n,b)$ of the table ${\mathcal R}$, with $n+b$ being automatically even, is $\langle \rho, \beta \rangle = \mathfrak{e}(\rho)[n,b]=\sum_{a,\, \partial a = 0} (F_{n})_{a,b} + \sum_{a,\, \partial a = 1} (F_{n-1})_{a,b}$, where the matrices $F_n$ are 
determined from the Chebyshev recurrence relation and the seed $F_1=\one$, $F_2=G$. In particular, $F_0=F_N=0$.
The scalar product between $\rho$ and the root $\alpha$, also localized at the position $(n,b)$ is $\langle \rho, \alpha \rangle = \mathfrak{f}(\rho)[n,b]$.
If $\alpha$ is positive it is enough to consider $\mathfrak{f^+}(\rho)$.
If the root is long, or if $G$ is simply-laced, we have $\mathfrak{f}(\rho)[n,b]=\mathfrak{e}(\rho)[n,b]$. Otherwise these two values differ by a scaling factor which is constant along a given orbit of the Coxeter element $c$ (a column of $\mathfrak{e}$): for a simple root $\alpha$ and corresponding coroot $\beta$, we have $\alpha = \delta\beta$ where the value of $\delta$ was recalled in section \ref{sec:nonADEessmat}, see also footnote \ref{foot: symmetrizingcoefficients}.
The classical (\ie $q=1$) Weyl denominator ${\mathcal D} = {\mathcal D}_{q=1}$ is   
 $\Pi_{{}^{n=1, b=1}_{n+b \, \text{even}}}^{n=N, b=r}\; \mathfrak{f^+}(\rho)[n,b]$.
If $G$ is simply-laced, this reads:
\[
{\mathcal D} =  \Pi_{{}^{n=1, b=1}_{n+b \, \text{even}}}^{n=N, \, b=r}\; \left( \sum_{a \, {\text with} \, \partial a=0} (F_{n})_{a,b} + \sum_{a \, {\text with} \, \partial a=1} (F_{n-1})_{a,b} \right)
\tag*{\ref{sec:weyldenominatorfromfusion}-1}
\]
If $G$ is not simply-laced, one has to multiply the above by a constant $\varpi$ equal to the product of known scaling factors: for all non-ADE Dynkin diagrams, the Coxeter number $N$ is even and the orbits of $c$ contain $N/2$ positive roots;
there is one orbit of short roots for $G_2$, two for $F_4$, one for $B_r$ and $r-1$ for $C_r$.
Using then $\varpi= \Pi_{s=1,\ldots r}\; \delta_s^{N/2}$ where $\delta_s$ as in \ref{sec:nonADEessmat} and footnote \ref{foot: symmetrizingcoefficients}, one obtains $\varpi = (1/3)^3$ for $G_2$, $(1/2)^{2 \times 6}$ for $F_4$, $(1/2)^r$ for $B_r$ and  $(1/2)^{r(r-1)}$ for $C_r$. We set $\varpi = 1$ for the ADE cases.

Permuting $F_n$ and  $F_{n-1}$ in the definition of tables $\mathfrak{e}(\omega)$ while imposing that $n+b$ should be odd rather than even 
amounts to change the orientation of the quiver of roots, \ie to use the bipartite Coxeter element $c^{-1}$ rather than $c$ to build the orbits.
For this reason permuting $F_n$ and  $F_{n-1}$ in the above expression while 
taking the  product over $n$ and $b$ with $n+b$ {\it odd} does not change the left hand side. In other words the list of elements of the $N \times r$ matrix $( \sum_{a=1}^r( (F_n)_{a,b} +  (F_{n-1})_{a,b}))$, indexed by $(n,b)$,
is obtained from the list of scalar products $\langle \rho, \beta \rangle$, where $\beta$ are positive coroots, by just doubling the multiplicity.
This remark allows us to remove all parity constraints in the previous expression for  ${\mathcal D}$, and obtain, for the square of the Weyl denominator, the identity:
\[
{\mathcal D}^2 \;= \varpi^2  \times \;\prod_{n=1}^{N} \prod_{b=1}^r \sum_{a=1}^r  \; (F_n+F_{n-1})_{a,b}
\tag*{\ref{sec:weyldenominatorfromfusion}-2}\label{bigprodeq}
\]

The quantum Weyl denominator ${\mathcal D}_q$ is  obtained without further ado by just replacing the value of each scalar product $\langle \rho, \alpha \rangle$ by its q-number analog $\langle \rho, \alpha \rangle_q$.
Warning:  in the non-ADE cases, one should be careful not to factor out the scaling factors $\delta_s$ -- see the example $F_4$ in the next section.

\section{The periodic quiver representation of the  Weyl vector}
\label{sec:periodicquiver}

The results of section \ref{sec: globaldim} require an explicit determination of the list of scalar products between the Weyl vector $\rho$ and all the roots (actually all the positive roots) of the chosen simple Lie group.
Of course, this can be obtained by a variety of means, and this is done in the appendix, using a direct approach, for the generic members of the series $A_r$, $B_r$, $C_r$, $D_r$. The same type of calculation can be done for the exceptional cases as well, but it is in the spirit of the paper to obtain the same information by making use of the correspondence between module-categories of type $SU(2)$ and ADE Dynkin diagrams (with the appropriate modifications needed to handle the non-ADE cases). 
In all cases, the determination of the matrices $F_n$, and therefore of the tables $\mathfrak{e}^{}_a$,  for all fundamental weights $\omega_a$, and of the tables ${}{\mathfrak{f}^{+}(\varrho)}$,  is almost immediate, given the adjacency matrix. We shall describe in some detail the example $E_6$, and to some extent the example $F_4$, in order to illustrate the discussion carried out in the last section, in particular in  \ref{sec:fusiontoscalarproducts}.
We shall then consider the other exceptional cases and give the expressions ${}{\mathfrak{f}^{+}(\varrho)}$ for the Weyl vector, from which the stated results for the Weyl denominator follow.
We shall also illustrate this technique when $G$ is not exceptional, by considering the cases $A_6, B_6, C_6, D_6$. For generic $r$, a proof of the fact that the classical  Weyl denominator is given by a Lie superfactorial could in principle also be obtained by evaluating the r.h.s. of equation \ref{bigprodeq} giving  ${\mathcal D}^2$. This is however not straightforward, as the calculation (see section\ref{sec:fusionmethodfornonADE}) involves non-trivial identifies between superfactorial functions, even in the relatively simple case of $A_r$ (for generic $r$).

 \bigskip
 
 In the non simply-laced cases the adjacency matrix $G$ is not symmetric and some of its elements are integers bigger than $1$. 
If $a$ and $b$ denote two vertices such that $G_{ab}\neq G_{ba}$ we draw one oriented edge from $a$ to $b$, with a multiplicity $G_{ab}$ displayed on the side of the chosen edge, and one oriented edge from $b$ to $a$, with its multiplicity $G_{ba}$. In the case $G_{ab}=G_{ba}$, we draw a single un-oriented edge between $a$ and $b$, and the multiplicity is omitted if it is equal to $1$. See the example of $F_4$ below.
\subsection{Two examples:  $E_6$ and $F_4$}
\subsubsection*{Case  $E_6$,
 \hskip 1.cm
\texorpdfstring{$G\,=\,$
{\tiny 
$\pmatrixcmd
{0 & 1 & 0 & 0 & 0 & 0 \\
 1 & 0 & 1 & 0 & 0 & 0 \\
 0 & 1 & 0 & 1 & 0 & 1 \\
 0 & 0 & 1 & 0 & 1 & 0 \\
 0 & 0 & 0 & 1 & 0 & 0 \\
 0 & 0 & 1 & 0 & 0 & 0}
$}
$ \Longrightarrow \mathfrak{e}^{+}(\varrho) =$
{\tiny 
$\pmatrixcmd
{ 1 &  \text{} & 1 &  \text{} & 1 &  \text{} \\
  \text{} & 3 &  \text{} & 3 &  \text{} & 2 \\
 2 &  \text{} & 7 &  \text{} & 2 &  \text{} \\
  \text{} & 6 &  \text{} & 6 &  \text{} & 5 \\
 4 &  \text{} & 10 &  \text{} & 4 &  \text{} \\
  \text{} & 8 &  \text{} & 8 &  \text{} & 5 \\
 4 &  \text{} & 11 &  \text{} & 4 &  \text{} \\
  \text{} & 7 &  \text{} & 7 &  \text{} & 6 \\
 3 &  \text{} & 9 &  \text{} & 3 &  \text{} \\
  \text{} & 5 &  \text{} & 5 &  \text{} & 3 \\
 2 &  \text{} & 4 &  \text{} & 2 &  \text{} \\
  \text{} & 1 &  \text{} & 1 &  \text{} & 1 \\}
$}
}{Mat ()}}


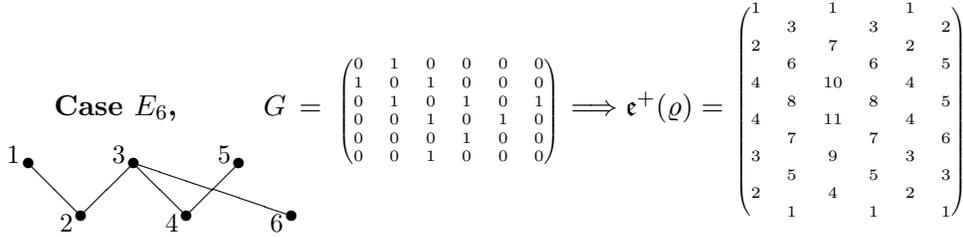
\unitlength 0.7mm
\begin{picture}(5,5)(25,0)
\put(20,20){\line(1,-1){10}}
\put(40,20){\line(-1,-1){10}}
\put(40,20){\line(1,-1){10}}
\put(60,20){\line(-1,-1){10}}

\put(40,20){\line(3,-1){30}}
\put(20,20){\circle*{2}} \put(16,20){1}
\put(40,20){\circle*{2}} \put(36,20){3}
\put(60,20){\circle*{2}} \put(56,20){5}

\put(30,10){\circle*{2}}\put(26,7){2}
\put(50,10){\circle*{2}}\put(46,7){4}
\put(70,10){\circle*{2}}\put(66,7){6}
\end{picture}


\label{Example:E6}

We call, as before, ${\omega_a}$ the fundamental weights, and ${\alpha_a}$ the simple roots. Their ordering $1,2,\ldots,6$ is specified by the nodes of the Dynkin diagram or by the adjacency matrix $G$. The fusion matrices $F_n$ are calculated from the Chebyshev recurrence formula, together with the seed $F_1=\one, F_2=G$. Here $N=N^\vee=12$ and the period is $2N=24$. The six $24 \times 6$ rectangular matrices $\mathfrak{e}_a$ (or the six $12 \times 6$ rectangular matrices $\mathfrak{e}_a^{+}$) are obtained as explained in section \ref{sec:simplylaced}.
Each fundamental weight  $\omega_a$ determines a table $\mathfrak{e}(\omega_a) = \mathfrak{e}_a$, \ie a function on the set ${\mathcal R}$; the $\mathfrak{e}_a^{+}$  are explicitly given below. 

\newenvironment{changemargin}[2]{\begin{list}{}{%
\setlength{\topsep}{0pt}%
\setlength{\leftmargin}{0pt}%
\setlength{\rightmargin}{0pt}%
\setlength{\listparindent}{\parindent}%
\setlength{\itemindent}{\parindent}%
\setlength{\parsep}{0pt plus 1pt}%
\addtolength{\leftmargin}{#1}%
\addtolength{\rightmargin}{#2}%
}\item }{\end{list}}

\begin{changemargin}{0.cm}{+2cm}
\begin{tabular}{cc}
\begin{tabular}{c}
\begin{tabular}{ccc}
{\tiny
$
\mathfrak{e}^{+}_1 =
\begin{array}{cccccc}
 1 & \text{} & 0 & \text{} & 0 & \text{} \\
 \text{} & 1 & \text{} & 0 & \text{} & 0 \\
 0 & \text{} & 1 & \text{} & 0 & \text{} \\
 \text{} & 0 & \text{} & 1 & \text{} & 1 \\
 0 & \text{} & 1 & \text{} & 1 & \text{} \\
 \text{} & 1 & \text{} & 1 & \text{} & 0 \\
 1 & \text{} & 1 & \text{} & 0 & \text{} \\
 \text{} & 1 & \text{} & 0 & \text{} & 1 \\
 0 & \text{} & 1 & \text{} & 0 & \text{} \\
 \text{} & 0 & \text{} & 1 & \text{} & 0 \\
 0 & \text{} & 0 & \text{} & 1 & \text{} \\
 \text{} & 0 & \text{} & 0 & \text{} & 0 \\
\end{array}
,
\quad
\mathfrak{e}^{+}_2 =
\begin{array}{cccccc}
 0 & \text{} & 0 & \text{} & 0 & \text{} \\
 \text{} & 1 & \text{} & 0 & \text{} & 0 \\
 1 & \text{} & 1 & \text{} & 0 & \text{} \\
 \text{} & 1 & \text{} & 1 & \text{} & 1 \\
 0 & \text{} & 2 & \text{} & 1 & \text{} \\
 \text{} & 1 & \text{} & 2 & \text{} & 1 \\
 1 & \text{} & 2 & \text{} & 1 & \text{} \\
 \text{} & 2 & \text{} & 1 & \text{} & 1 \\
 1 & \text{} & 2 & \text{} & 0 & \text{} \\
 \text{} & 1 & \text{} & 1 & \text{} & 1 \\
 0 & \text{} & 1 & \text{} & 1 & \text{} \\
 \text{} & 0 & \text{} & 1 & \text{} & 0 \\
\end{array}
,
\quad
\mathfrak{e}^{+}_3 =
\begin{array}{cccccc}
 0 & \text{} & 1 & \text{} & 0 & \text{} \\
 \text{} & 1 & \text{} & 1 & \text{} & 1 \\
 1 & \text{} & 2 & \text{} & 1 & \text{} \\
 \text{} & 2 & \text{} & 2 & \text{} & 1 \\
 1 & \text{} & 3 & \text{} & 1 & \text{} \\
 \text{} & 2 & \text{} & 2 & \text{} & 2 \\
 1 & \text{} & 3 & \text{} & 1 & \text{} \\
 \text{} & 2 & \text{} & 2 & \text{} & 1 \\
 1 & \text{} & 2 & \text{} & 1 & \text{} \\
 \text{} & 1 & \text{} & 1 & \text{} & 1 \\
 0 & \text{} & 1 & \text{} & 0 & \text{} \\
 \text{} & 0 & \text{} & 0 & \text{} & 0 \\
\end{array}
$
}
\end{tabular}
\\
\begin{tabular}{ccc}
 &  &   \\
\end{tabular}
\\
\begin{tabular}{ccc}
{\tiny
$
\mathfrak{e}^{+}_4 =
\begin{array}{cccccc}
 0 & \text{} & 0 & \text{} & 0 & \text{} \\
 \text{} & 0 & \text{} & 1 & \text{} & 0 \\
 0 & \text{} & 1 & \text{} & 1 & \text{} \\
 \text{} & 1 & \text{} & 1 & \text{} & 1 \\
 1 & \text{} & 2 & \text{} & 0 & \text{} \\
 \text{} & 2 & \text{} & 1 & \text{} & 1 \\
 1 & \text{} & 2 & \text{} & 1 & \text{} \\
 \text{} & 1 & \text{} & 2 & \text{} & 1 \\
 0 & \text{} & 2 & \text{} & 1 & \text{} \\
 \text{} & 1 & \text{} & 1 & \text{} & 1 \\
 1 & \text{} & 1 & \text{} & 0 & \text{} \\
 \text{} & 1 & \text{} & 0 & \text{} & 0 \\
\end{array}
,
\quad
\mathfrak{e}^{+}_5 =
\begin{array}{cccccc}
 0 & \text{} & 0 & \text{} & 1 & \text{} \\
 \text{} & 0 & \text{} & 1 & \text{} & 0 \\
 0 & \text{} & 1 & \text{} & 0 & \text{} \\
 \text{} & 1 & \text{} & 0 & \text{} & 1 \\
 1 & \text{} & 1 & \text{} & 0 & \text{} \\
 \text{} & 1 & \text{} & 1 & \text{} & 0 \\
 0 & \text{} & 1 & \text{} & 1 & \text{} \\
 \text{} & 0 & \text{} & 1 & \text{} & 1 \\
 0 & \text{} & 1 & \text{} & 0 & \text{} \\
 \text{} & 1 & \text{} & 0 & \text{} & 0 \\
 1 & \text{} & 0 & \text{} & 0 & \text{} \\
 \text{} & 0 & \text{} & 0 & \text{} & 0 \\
\end{array}
,
\quad
\mathfrak{e}^{+}_6 =
\begin{array}{cccccc}
 0 & \text{} & 0 & \text{} & 0 & \text{} \\
 \text{} & 0 & \text{} & 0 & \text{} & 1 \\
 0 & \text{} & 1 & \text{} & 0 & \text{} \\
 \text{} & 1 & \text{} & 1 & \text{} & 0 \\
 1 & \text{} & 1 & \text{} & 1 & \text{} \\
 \text{} & 1 & \text{} & 1 & \text{} & 1 \\
 0 & \text{} & 2 & \text{} & 0 & \text{} \\
 \text{} & 1 & \text{} & 1 & \text{} & 1 \\
 1 & \text{} & 1 & \text{} & 1 & \text{} \\
 \text{} & 1 & \text{} & 1 & \text{} & 0 \\
 0 & \text{} & 1 & \text{} & 0 & \text{} \\
 \text{} & 0 & \text{} & 0 & \text{} & 1 \\
\end{array}
$
}
\end{tabular}
\end{tabular}
&
\begin{tabular}{c}
$
\begin{array}{c}
 {}  
\end{array}
$
\\
{}
\\
$
\begin{array}{c}
 {}  
\end{array}
$
\end{tabular}
\end{tabular}
\end{changemargin}

For illustration consider the simple root $\alpha_3$. From the third line of $G$ (or from the Cartan matrix $2\one - G$) we see that 
$\alpha_3 = 2 \omega_3- \omega_2  - \omega_4 - \omega_6$; the table $\mathfrak{e}(\alpha_3)= - 2 \mathfrak{e}_3 - \mathfrak{e}_2 -\mathfrak{e}_4 -\mathfrak{e}_6$ is displayed on figure \ref{rhoofalpha3forE6}.
Notice that it has a coefficient $+2$ on position $(n=1, a=3)$; this specifies the location of that particular root on ${\mathcal R}$. In the same way, one can determine the position, on $\mathcal R$,  of the $72$ roots of $E_6$, and displays the results on a $24 \times 6$ table (see figure \ref{table:rootsofE6}).
In this example, even and odd simple roots appear respectively on lines $n=1$ and $n=12$ of the table.
\begin{figure}[htpb]
\begin{minipage}{0.4\linewidth}
\centering
\begin{tikzpicture}
\tikzstyle{every node}=[black]
\matrix [matrix of math nodes, row sep=3.mm](Mat)
{
|[minimum width=3em]|  & |[minimum width=3em]| & |[minimum width=3em]| & |[minimum width=3em]|  & |[minimum width=3em]|    & |[minimum width=3em]| \\
 0 & & {\bf 2} & & 0 & \\
 & 1 & & 1 & & 1 \\
 1 & & 1 & & 1 & \\
 & 1 & & 1 & & 0 \\
 0 & & 1 & & 0 & \\
 & 0 & & 0 & & 1 \\
 0 & & 0 & & 0 & \\
 & 0 & & 0 & & -1 \\
 0 & & -1 & & 0 & \\
 & -1 & & -1 & & 0 \\
 -1 & & -1 & & -1 & \\
 & -1 & & -1 & & -1 \\
 0 & & -2 & & 0 & \\
 & -1 & & -1 & & -1 \\
 -1 & & -1 & & -1 & \\
 & -1 & & -1 & & 0 \\
 0 & & -1 & & 0 & \\
 & 0 & & 0 & & -1 \\
 0 & & 0 & & 0 & \\
 & 0 & & 0 & & 1 \\
 0 & & 1 & & 0 & \\
 & 1 & & 1 & & 0 \\
 1 & & 1 & & 1 & \\
 & 1 & & 1 & & 1 \\
};
\draw[very thick, red] (Mat-2-1) -- (Mat-3-2)node[]{};
\draw[very thick, red] (Mat-3-2) -- (Mat-2-3)node[]{};
\draw[very thick, red] (Mat-2-3) -- (Mat-3-4)node[]{};
\draw[very thick, red] (Mat-3-4) -- (Mat-2-5)node[]{};
\draw[very thick, red] (Mat-2-3) -- (Mat-3-6)node[]{};
\draw (Mat-4-1) -- (Mat-3-2)node[]{};
\draw (Mat-3-2) -- (Mat-4-3)node[]{};
\draw (Mat-4-3) -- (Mat-3-4)node[]{};
\draw (Mat-3-4) -- (Mat-4-5)node[]{};
\draw (Mat-4-3) -- (Mat-3-6)node[]{};
\draw (Mat-4-1) -- (Mat-5-2)node[]{};
\draw (Mat-5-2) -- (Mat-4-3)node[]{};
\draw (Mat-4-3) -- (Mat-5-4)node[]{};
\draw (Mat-5-4) -- (Mat-4-5)node[]{};
\draw (Mat-4-3) -- (Mat-5-6)node[]{};
\draw (Mat-6-1) -- (Mat-5-2)node[]{};
\draw (Mat-5-2) -- (Mat-6-3)node[]{};
\draw (Mat-6-3) -- (Mat-5-4)node[]{};
\draw (Mat-5-4) -- (Mat-6-5)node[]{};
\draw (Mat-6-3) -- (Mat-5-6)node[]{};
\end{tikzpicture}
\caption{$\mathfrak{e}(\alpha_3)$}
\label{rhoofalpha3forE6}
\end{minipage}
\hspace{3.0cm}
\noindent\rotatebox{-90}{%
\begin{minipage}{5.cm}%
\begin{changemargin}{-10.0cm}{+3.0cm}
{\tiny
$
{\setlength{\arraycolsep}{2pt}
\begin{array}{cccccc}
 \alpha _1 &  \text{} & \alpha _3 &  \text{} & \alpha _5 &  \text{} \\
  \text{} & \alpha _1+\alpha _2+\alpha _3 &  \text{} & \alpha _3+\alpha _4+\alpha _5 &  \text{} & \alpha _3+\alpha _6 \\
 \alpha _2+\alpha _3 &  \text{} & \alpha _1+\alpha _2+2 \alpha _3+\alpha _4+\alpha _5+\alpha _6 &  \text{} & \alpha _3+\alpha _4 &  \text{} \\
  \text{} & \alpha _2+2 \alpha _3+\alpha _4+\alpha _5+\alpha _6 &  \text{} & \alpha _1+\alpha _2+2 \alpha _3+\alpha _4+\alpha _6 &  \text{} & \alpha _1+\alpha _2+\alpha _3+\alpha _4+\alpha _5 \\
 \alpha _3+\alpha _4+\alpha _5+\alpha _6 &  \text{} & \alpha _1+2 \alpha _2+3 \alpha _3+2 \alpha _4+\alpha _5+\alpha _6 &  \text{} & \alpha _1+\alpha _2+\alpha _3+\alpha _6 &  \text{} \\
  \text{} & \alpha _1+\alpha _2+2 \alpha _3+2 \alpha _4+\alpha _5+\alpha _6 &  \text{} & \alpha _1+2 \alpha _2+2 \alpha _3+\alpha _4+\alpha _5+\alpha _6 &  \text{} & \alpha _2+2 \alpha _3+\alpha
   _4+\alpha _6 \\
 \alpha _1+\alpha _2+\alpha _3+\alpha _4 &  \text{} & \alpha _1+2 \alpha _2+3 \alpha _3+2 \alpha _4+\alpha _5+2 \alpha _6 &  \text{} & \alpha _2+\alpha _3+\alpha _4+\alpha _5 &  \text{} \\
  \text{} & \alpha _1+2 \alpha _2+2 \alpha _3+\alpha _4+\alpha _6 &  \text{} & \alpha _2+2 \alpha _3+2 \alpha _4+\alpha _5+\alpha _6 &  \text{} & \alpha _1+\alpha _2+\alpha _3+\alpha _4+\alpha _5+\alpha
   _6 \\
 \alpha _2+\alpha _3+\alpha _6 &  \text{} & \alpha _1+2 \alpha _2+2 \alpha _3+2 \alpha _4+\alpha _5+\alpha _6 &  \text{} & \alpha _3+\alpha _4+\alpha _6 &  \text{} \\
  \text{} & \alpha _2+\alpha _3+\alpha _4+\alpha _5+\alpha _6 &  \text{} & \alpha _1+\alpha _2+\alpha _3+\alpha _4+\alpha _6 &  \text{} & \alpha _2+\alpha _3+\alpha _4 \\
 \alpha _4+\alpha _5 &  \text{} & \alpha _2+\alpha _3+\alpha _4+\alpha _6 &  \text{} & \alpha _1+\alpha _2 &  \text{} \\
  \text{} & \alpha _4 &  \text{} & \alpha _2 &  \text{} & \alpha _6 \\
 -\alpha _5 &  \text{} & -\alpha _3 &  \text{} & -\alpha _1 &  \text{} \\
  \text{} & -\alpha _3-\alpha _4-\alpha _5 &  \text{} & -\alpha _1-\alpha _2-\alpha _3 &  \text{} & -\alpha _3-\alpha _6 \\
 -\alpha _3-\alpha _4 &  \text{} & -\alpha _1-\alpha _2-2 \alpha _3-\alpha _4-\alpha _5-\alpha _6 &  \text{} & -\alpha _2-\alpha _3 &  \text{} \\
  \text{} & -\alpha _1-\alpha _2-2 \alpha _3-\alpha _4-\alpha _6 &  \text{} & -\alpha _2-2 \alpha _3-\alpha _4-\alpha _5-\alpha _6 &  \text{} & -\alpha _1-\alpha _2-\alpha _3-\alpha _4-\alpha _5 \\
 -\alpha _1-\alpha _2-\alpha _3-\alpha _6 &  \text{} & -\alpha _1-2 \alpha _2-3 \alpha _3-2 \alpha _4-\alpha _5-\alpha _6 &  \text{} & -\alpha _3-\alpha _4-\alpha _5-\alpha _6 &  \text{} \\
  \text{} & -\alpha _1-2 \alpha _2-2 \alpha _3-\alpha _4-\alpha _5-\alpha _6 &  \text{} & -\alpha _1-\alpha _2-2 \alpha _3-2 \alpha _4-\alpha _5-\alpha _6 &  \text{} & -\alpha _2-2 \alpha _3-\alpha
   _4-\alpha _6 \\
 -\alpha _2-\alpha _3-\alpha _4-\alpha _5 &  \text{} & -\alpha _1-2 \alpha _2-3 \alpha _3-2 \alpha _4-\alpha _5-2 \alpha _6 &  \text{} & -\alpha _1-\alpha _2-\alpha _3-\alpha _4 &  \text{} \\
  \text{} & -\alpha _2-2 \alpha _3-2 \alpha _4-\alpha _5-\alpha _6 &  \text{} & -\alpha _1-2 \alpha _2-2 \alpha _3-\alpha _4-\alpha _6 &  \text{} & -\alpha _1-\alpha _2-\alpha _3-\alpha _4-\alpha
   _5-\alpha _6 \\
 -\alpha _3-\alpha _4-\alpha _6 &  \text{} & -\alpha _1-2 \alpha _2-2 \alpha _3-2 \alpha _4-\alpha _5-\alpha _6 &  \text{} & -\alpha _2-\alpha _3-\alpha _6 &  \text{} \\
  \text{} & -\alpha _1-\alpha _2-\alpha _3-\alpha _4-\alpha _6 &  \text{} & -\alpha _2-\alpha _3-\alpha _4-\alpha _5-\alpha _6 &  \text{} & -\alpha _2-\alpha _3-\alpha _4 \\
 -\alpha _1-\alpha _2 &  \text{} & -\alpha _2-\alpha _3-\alpha _4-\alpha _6 &  \text{} & -\alpha _4-\alpha _5 &  \text{} \\
  \text{} & -\alpha _2 &  \text{} & -\alpha _4 &  \text{} & -\alpha _6 \\
\end{array} 
}%
$
}%
\caption{\label{table:rootsofE6} The $E_6$ quiver of roots}
\end{changemargin}
\end{minipage}
}%
\end{figure}
The scalar products between an arbitrary weight $\omega = \sum_a \lambda_a \omega_a$ and the roots (resp. positive roots) are given by the entries of the table $\mathfrak{e}(\omega) = \sum_a \lambda_a \, \mathfrak{e}_a$ (resp. of the table $\mathfrak{e}^{+}(\omega)$).
For instance, the position of the non-simple root $\alpha=\alpha_2 + \alpha_3 + \alpha_4 + \alpha_6$ is $(n=11,\, a=3)$.  Its scalar products with the fundamental weights $\omega_a$, for $a=1,2,\ldots,6$, are given by the numbers sitting in position $(n=11,\, a=3)$ of the tables $\omega_a$, and its scalar product with the root $\alpha_3$, for instance, is $-1$, as we see from the table $\mathfrak{e}(\alpha_3)$.

It is not necessary for the purpose of this paper to determine the position of all the roots on the set $\mathcal R$, not even the position of the simple roots\footnote{
Calling ${\omega_a}$ the basis of fundamental weights and  ${\beta_b}$ the dual basis of simple coroots with respect to the scalar product $\langle \, \cdot , \cdot \, \rangle$,  \ie $\langle \omega_a , \beta_b \rangle = \delta_{ab}$, we have $\varrho = \sum_a \, \omega_a$, so that  $\langle \varrho , \beta_a \rangle =1$ for all $a \in G$. This determines the positions of simple coroots (among all positive coroots) on the vertices of the half-quiver $\mathcal R^{+}$. The reader can check this property by comparing the entries $1$ in ${\mathfrak e}^{+}(\varrho)$ with the position of the simple coroots $\beta_a$ in the table giving the positions of all coroots. Of course, for ADE diagrams, roots and coroots coincide ($\alpha_a=\beta_a$).}.
This would be needed if we wanted to identify the {\it individual} numbers appearing in the arrays $\mathfrak{e}(\omega)$ as specific scalar products.
We gave this information -- in the case of $E_6$ -- for illustration  but what we only need is the list of entries of the table $\mathfrak{e}^{+}(\varrho) =  \sum_a \, \mathfrak{e}^{+}_a$, where $\varrho = \sum_a \, \omega_a$ is the Weyl vector.  
 This table  can be readily calculated from the fusion matrices, as explained, and it is displayed next to the Dynkin diagram of $E_6$.  
  By multiplying together all its entries, one obtains the Weyl denominator  ${\mathcal D} = 1!  \, 4!  \, 5! \,  7!  \,  8!  \,  11! = \prod_{s \in {\mathcal Exp(E_6)}} \,  s!  = \mathrm{\sfac}_{E_6}$, as it was announced.
  Since $\langle \omega_a , \beta_b \rangle = \delta_{ab}$, one can check that $\mathfrak{e}^{}(\varrho)$ can also be obtained from table \ref{table:rootsofE6} by replacing all the simple roots $\alpha_a$ by $1$.
 
One should think of $\mathcal R$ as an oriented graph: its vertices have positive integer levels $n$ counted from the top ($n=1$); its edges, that are oriented from top to bottom, are drawn between vertices of ${\mathcal R}$ having levels differing by one and relate pairs of vertices  that correspond to neighbors on the Dynkin diagram $G$ (neighbors on $G$ are automatically of opposite parity for the given $\ZZ_2$ grading), see also the remark made at the end of section \ref{sec:fusiontoscalarproducts}.
We only display the edges relating the first five levels of the table  $\mathfrak{e}(\alpha_3)$:  they follow the pattern determined by the folded Dynkin diagram drawn next to the adjacency matrix;  the pattern is then reflected horizontally every two steps when going downwards, and the whole structure is periodic, the last level being connected with the first.

\subsubsection*{Case  $F_4$, 
 \hskip 0.2cm
\texorpdfstring{$G\,=\,$
{\tiny 
$\pmatrixcmd
{
 0 & 1 & 0 & 0 \\
 1 & 0 & 2 & 0 \\
 0 & 1 & 0 & 1 \\
 0 & 0 & 1 & 0 \\
}
$}
$ \Rightarrow \mathfrak{e}^{+}(\varrho) =$
{\tiny 
$\pmatrixcmd
{ 1 & \text{} & 1 & \text{} \\
 \text{} & 3 & \text{} & 2 \\
 2 & \text{} & 7 & \text{} \\
 \text{} & 6 & \text{} & 5 \\
 4 & \text{} & 10 & \text{} \\
 \text{} & 8 & \text{} & 5 \\
 4 & \text{} & 11 & \text{} \\
 \text{} & 7 & \text{} & 6 \\
 3 & \text{} & 9 & \text{} \\
 \text{} & 5 & \text{} & 3 \\
 2 & \text{} & 4 & \text{} \\
 \text{} & 1 & \text{} & 1 }
$}
$ \Rightarrow {}{\mathfrak{f}^{+}(\varrho)} =$
{\tiny 
$\pmatrixcmd
{ 1 & \text{} & 1/2 & \text{} \\
 \text{} & 3 & \text{} & 2/2 \\
 2 & \text{} & 7/2 & \text{} \\
 \text{} & 6 & \text{} & 5/2 \\
 4 & \text{} & 10/2 & \text{} \\
 \text{} & 8 & \text{} & 5/2 \\
 4 & \text{} & 11/2 & \text{} \\
 \text{} & 7 & \text{} & 6/2 \\
 3 & \text{} & 9/2 & \text{} \\
 \text{} & 5 & \text{} & 3/2 \\
 2 & \text{} & 4/2 & \text{} \\
 \text{} & 1 & \text{} & 1/2 }
$}
 }{Mat ()}}


\unitlength 0.7mm
\begin{picture}(5,5)(25,0)
\put(20,20){\line(1,-1){10}}
\put(40,20){\line(-1,-1){10}}
\put(40,20){\line(1,-1){10}}

\qbezier(40, 20)(40,10)(30,10)

\put(30,18){$\nearrow^2$}
\put(35,7){$\swarrow^1$}
\put(20,20){\circle*{2}}
\put(40,20){\circle*{2}}

\put(30,10){\circle*{2}}
\put(50,10){\circle*{2}}
\end{picture}

The Chebyshev recurrence relation using the adjacency matrix $G$ of $F_4$ has a period $2 N$, with $N = 12$ (here $N^\vee=9$), so,  from the definition of matrices $\mathfrak{e}^{+}$, we obtain

{\tiny
$$
\mathfrak{e}^{+}_1 =
\begin{array}{cccc}
 1 & \text{} & 0 & \text{} \\
 \text{} & 1 & \text{} &0 \\
 0 & \text{} & 2 & \text{} \\
 \text{} & 1 & \text{} & 2 \\
 1 & \text{} & 2 & \text{} \\
 \text{} & 2 & \text{} & 0\\
 1 & \text{} & 2 & \text{} \\
 \text{} & 1 & \text{} & 2 \\
 0 & \text{} & 2 & \text{} \\
 \text{} & 1 & \text{} & 0 \\
 1 & \text{} & 0 & \text{} \\
 \text{} & 0 & \text{} & 0 
\end{array}
,
\quad
\mathfrak{e}^{+}_2 =
\begin{array}{cccc}
 0 & \text{} & 0 & \text{} \\
 \text{} & 1 & \text{} & 0 \\
 1 & \text{} & 2 & \text{} \\
 \text{} & 2 & \text{} & 2 \\
 1 & \text{} & 4 & \text{} \\
 \text{} & 3 & \text{} & 2 \\
 2 & \text{} & 4 & \text{} \\
 \text{} & 3 & \text{} & 2 \\
 1 & \text{} & 4 & \text{} \\
 \text{} & 2 & \text{} & 2 \\
 1 & \text{} & 2 & \text{} \\
 \text{} & 1 & \text{} & 0 
\end{array}
,
\quad
\mathfrak{e}^{+}_3 =
\begin{array}{cccc}
 0 & \text{} & 1 & \text{} \\
 \text{} & 1 & \text{} & 1 \\
 1 & \text{} & 2 & \text{} \\
 \text{} & 2 & \text{} & 1 \\
 1 & \text{} & 3 & \text{} \\
 \text{} & 2 & \text{} & 2 \\
 1 & \text{} & 3 & \text{} \\
 \text{} & 2 & \text{} & 1 \\
 1 & \text{} & 2 & \text{} \\
 \text{} & 1 & \text{} & 1 \\
 0 & \text{} & 1 & \text{} \\
 \text{} & 0 & \text{} & 0 
\end{array}
,
\quad
\mathfrak{e}^{+}_4 =
\begin{array}{cccc}
 0 & \text{} & 0 & \text{} \\
 \text{} & 0 & \text{} & 1 \\
 0 & \text{} & 1 & \text{} \\
 \text{} & 1 & \text{} & 0 \\
 1 & \text{} & 1 & \text{} \\
 \text{} & 1 & \text{} & 1 \\
 0 & \text{} & 2 & \text{} \\
 \text{} & 1 & \text{} & 1 \\
 1 & \text{} & 1 & \text{} \\
 \text{} & 1 & \text{} & 0 \\
 0 & \text{} & 1 & \text{} \\
 \text{} & 0 & \text{} & 1 \\
\end{array}
$$
}

The scalar products between the Weyl vector  $\varrho = \sum_a \, \omega_a$ and the positive coroots are obtained as $\mathfrak{e}^{+}(\varrho) =  \sum_a \, \mathfrak{e}^{+}_a$, and displayed above.
The scalar products between the Weyl vector and the  positive roots are obtained by scaling the last two columns of $\mathfrak{e}^{+}(\varrho)$  with a factor $1/2$. 
By multiplying together all values of ${}{\mathfrak{f}^{+}(\varrho)}$, one obtains ${\mathcal D} = 1!  \, 5! \,  7!  \,  11!/2^{12} =  (1/2^{12}) \prod_{s \in {\mathcal Exp(F_4)}} \,  s!  = \mathrm{\sfac}_{F_4}$.
Notice that the quantum Weyl denominator is obtained from the classical one by replacing the scalar products {\sl themselves} by their q-analogs, so that for instance, $7/2$ becomes $[7/2]_q$, and this is not equal to  $[7]_q/[2]_q$. By multiplying the q-analogs of all entries of the table ${}{\mathfrak{f}^{+}(\varrho)}$,  we recover the expression of the quantum superfactorial function of type $F_4$ that was given in section \ref{sec:globaldimforfusion}.

\subsection{The other exceptional cases: $G_2$, $E_7$, $E_8$}

The periodic quiver representation of the Weyl vector \ie the table displaying $\mathfrak{e}^{}(\varrho)$ when $G$ is simply-laced,  and ${}{\mathfrak{f}^{}(\varrho)}$ when it is not,
comes immediately from the tables  $\mathfrak{e}_a$ obtained from the fusion matrices $F_n$.
Here, we are only interested in the half-quiver of positive roots: 
the results for the exceptional cases $G_2$, $E_7$ and $E_8$ are gathered in table \ref{tab:rhoG2E7E8}  (the cases $E_6$ and $F_4$ have already been described). 
The corresponding values for the classical and quantum Weyl denominators are then obtained by multiplying together the entries of $\mathfrak{f}^{+}(\varrho)$,
they give the classical or quantum Lie superfactorials defined in sections \ref{sec: mainresult} and \ref{sec:globaldimforfusion}.

\begin{table}
{\tiny
$$
\begin{array}{ccc}

G_2
\begin{array}{c}

\unitlength 0.7mm
\begin{picture}(5,5)(20,-5)
\put(20,10){\line(1,-1){10}}
\qbezier(20, 10)(40,10)(30, 0)
\put(12,2){$\nwarrow^3$}
\put(35,7){$\searrow^1$}
\put(30,0){\circle*{2}}

\put(20,10){\circle*{2}}
\end{picture}

\\

\begin{array}{cc}
 1 & \text{} \\
 \text{} & \frac{4}{3} \\
 3 & \text{} \\
 \text{} & \frac{5}{3} \\
 2 & \text{} \\
 \text{} & \frac{1}{3} \\
\end{array}

\end{array}

&

E_7
\begin{array}{c}

\unitlength 0.7mm
\begin{picture}(5,5)(50,0)
\put(20,20){\line(1,-1){10}}
\put(40,20){\line(-1,-1){10}}
\put(40,20){\line(1,-1){10}}
\put(60,20){\line(-1,-1){10}}
\put(60,20){\line(1,-1){10}}

\put(50,10){\line(3,1){30}}
\put(20,20){\circle*{2}}
\put(40,20){\circle*{2}}
\put(60,20){\circle*{2}}
\put(80,20){\circle*{2}}

\put(30,10){\circle*{2}}
\put(50,10){\circle*{2}}
\put(70,10){\circle*{2}}
\end{picture}

\\

\begin{array}{ccccccc}
 1 & \text{} & 1 & \text{} & 1 & \text{} & 1 \\
 \text{} & 3 & \text{} & 4 & \text{} & 2 & \text{} \\
 2 & \text{} & 6 & \text{} & 5 & \text{} & 3 \\
 \text{} & 5 & \text{} & 10 & \text{} & 3 & \text{} \\
 3 & \text{} & 9 & \text{} & 8 & \text{} & 7 \\
 \text{} & 7 & \text{} & 14 & \text{} & 5 & \text{} \\
 4 & \text{} & 12 & \text{} & 11 & \text{} & 7 \\
 \text{} & 9 & \text{} & 16 & \text{} & 6 & \text{} \\
 5 & \text{} & 13 & \text{} & 11 & \text{} & 9 \\
 \text{} & 9 & \text{} & 17 & \text{} & 5 & \text{} \\
 4 & \text{} & 13 & \text{} & 11 & \text{} & 8 \\
 \text{} & 8 & \text{} & 15 & \text{} & 6 & \text{} \\
 4 & \text{} & 10 & \text{} & 10 & \text{} & 7 \\
 \text{} & 6 & \text{} & 12 & \text{} & 4 & \text{} \\
 2 & \text{} & 8 & \text{} & 6 & \text{} & 5 \\
 \text{} & 4 & \text{} & 7 & \text{} & 2 & \text{} \\
 2 & \text{} & 3 & \text{} & 3 & \text{} & 2 \\
 \text{} & 1 & \text{} & 1 & \text{} & 1 & \text{} \\
\end{array}

\end{array}

&

E_8
\begin{array}{c}

\unitlength 0.7mm
\begin{picture}(5,5)(50,0)
\put(20,20){\line(1,-1){10}}
\put(40,20){\line(-1,-1){10}}
\put(40,20){\line(1,-1){10}}
\put(60,20){\line(-1,-1){10}}
\put(60,20){\line(1,-1){10}}
\put(80,20){\line(-1,-1){10}}

\put(60,20){\line(3,-1){30}}
\put(20,20){\circle*{2}}
\put(40,20){\circle*{2}}
\put(60,20){\circle*{2}}
\put(80,20){\circle*{2}}

\put(30,10){\circle*{2}}
\put(50,10){\circle*{2}}
\put(70,10){\circle*{2}}
\put(90,10){\circle*{2}}
\end{picture}
\\
\begin{array}{cccccccc}
 1 & \text{} & 1 & \text{} & 1 & \text{} & 1 & \text{} \\
 \text{} & 3 & \text{} & 3 & \text{} & 3 & \text{} & 2 \\
 2 & \text{} & 5 & \text{} & 7 & \text{} & 2 & \text{} \\
 \text{} & 4 & \text{} & 9 & \text{} & 6 & \text{} & 5 \\
 2 & \text{} & 8 & \text{} & 13 & \text{} & 4 & \text{} \\
 \text{} & 6 & \text{} & 12 & \text{} & 11 & \text{} & 8 \\
 4 & \text{} & 10 & \text{} & 18 & \text{} & 7 & \text{} \\
 \text{} & 8 & \text{} & 16 & \text{} & 14 & \text{} & 10 \\
 4 & \text{} & 14 & \text{} & 22 & \text{} & 7 & \text{} \\
 \text{} & 10 & \text{} & 20 & \text{} & 15 & \text{} & 12 \\
 6 & \text{} & 16 & \text{} & 25 & \text{} & 8 & \text{} \\
 \text{} & 12 & \text{} & 21 & \text{} & 18 & \text{} & 13 \\
 6 & \text{} & 17 & \text{} & 27 & \text{} & 10 & \text{} \\
 \text{} & 11 & \text{} & 23 & \text{} & 19 & \text{} & 14 \\
 5 & \text{} & 17 & \text{} & 29 & \text{} & 9 & \text{} \\
 \text{} & 11 & \text{} & 23 & \text{} & 19 & \text{} & 15 \\
 6 & \text{} & 17 & \text{} & 28 & \text{} & 10 & \text{} \\
 \text{} & 12 & \text{} & 22 & \text{} & 19 & \text{} & 13 \\
 6 & \text{} & 17 & \text{} & 26 & \text{} & 9 & \text{} \\
 \text{} & 11 & \text{} & 21 & \text{} & 16 & \text{} & 13 \\
 5 & \text{} & 15 & \text{} & 24 & \text{} & 7 & \text{} \\
 \text{} & 9 & \text{} & 18 & \text{} & 15 & \text{} & 11 \\
 4 & \text{} & 12 & \text{} & 20 & \text{} & 8 & \text{} \\
 \text{} & 7 & \text{} & 14 & \text{} & 13 & \text{} & 9 \\
 3 & \text{} & 9 & \text{} & 16 & \text{} & 5 & \text{} \\
 \text{} & 5 & \text{} & 11 & \text{} & 8 & \text{} & 7 \\
 2 & \text{} & 7 & \text{} & 10 & \text{} & 3 & \text{} \\
 \text{} & 4 & \text{} & 6 & \text{} & 5 & \text{} & 3 \\
 2 & \text{} & 3 & \text{} & 4 & \text{} & 2 & \text{} \\
 \text{} & 1 & \text{} & 1 & \text{} & 1 & \text{} & 1 \\
\end{array}
\end{array}
\end{array}
$$
}
\caption{Values of the Weyl vector on the half-quiver of positive roots. Cases $G_2$, $E_7$, $E_8$.}
\label{tab:rhoG2E7E8}
\end{table}
 
\subsection{Non-exceptional cases: the examples $A_6$, $B_6$, $C_6$, $D_6$}
\label{sec:fusionmethodfornonADE}

Like for the exceptional cases, for every specific choice of the rank $r$,  obtaining the matrices $F_n$ and the  tables $\mathfrak{f}_a$ for the classical series is a straightforward operation. 
Again, one can use these tables to express the corresponding Weyl denominator as a Lie superfactorial.
The periodic quiver representation of the Weyl vector is displayed below,  for the cases $A_6$, $B_6$, $C_6$, $D_6$,  in table \ref{tab:rhoA6B6C6D6}. 

\begin{table}
\vspace*{1. cm}
{\tiny
$$
\begin{array}{cccc}
A_6
\begin{array}{c}
\unitlength 0.5mm
\begin{picture}(5,5)(45,0)
\put(20,20){\line(1,-1){10}}
\put(40,20){\line(-1,-1){10}}
\put(40,20){\line(1,-1){10}}
\put(60,20){\line(-1,-1){10}}
\put(60,20){\line(1,-1){10}}
\put(20,20){\circle*{2}}
\put(40,20){\circle*{2}}
\put(60,20){\circle*{2}}

\put(30,10){\circle*{2}}
\put(50,10){\circle*{2}}
\put(70,10){\circle*{2}}
\end{picture}

\\

\begin{array}{cccccc}
 1 & \text{} & 1 & \text{} & 1 & \text{} \\
 \text{} & 3 & \text{} & 3 & \text{} & 2 \\
 2 & \text{} & 5 & \text{} & 4 & \text{} \\
 \text{} & 4 & \text{} & 6 & \text{} & 2 \\
 2 & \text{} & 5 & \text{} & 4 & \text{} \\
 \text{} & 3 & \text{} & 3 & \text{} & 2 \\
 1 & \text{} & 1 & \text{} & 1 & \text{} \\
\end{array}

\end{array}

&

D_6
\begin{array}{c}

\unitlength 0.5mm
\begin{picture}(5,5)(45,0)
\put(20,20){\line(1,-1){10}}
\put(40,20){\line(-1,-1){10}}
\put(40,20){\line(1,-1){10}}
\put(60,20){\line(-1,-1){10}}

\put(50,10){\line(3,1){30}}
\put(20,20){\circle*{2}}
\put(40,20){\circle*{2}}
\put(60,20){\circle*{2}}
\put(80,20){\circle*{2}}

\put(30,10){\circle*{2}}
\put(50,10){\circle*{2}}
\end{picture}

\\

\begin{array}{cccccc}
 1 & \text{} & 1 & \text{} & 1 & 1 \\
 \text{} & 3 & \text{} & 4 & \text{} & \text{} \\
 2 & \text{} & 6 & \text{} & 3 & 3 \\
 \text{} & 5 & \text{} & 8 & \text{} & \text{} \\
 3 & \text{} & 7 & \text{} & 5 & 5 \\
 \text{} & 5 & \text{} & 9 & \text{} & \text{} \\
 2 & \text{} & 7 & \text{} & 4 & 4 \\
 \text{} & 4 & \text{} & 6 & \text{} & \text{} \\
 2 & \text{} & 3 & \text{} & 2 & 2 \\
 \text{} & 1 & \text{} & 1 & \text{} & \text{} \\
\end{array}

\end{array}

&

B_6
\begin{array}{c}

\unitlength 0.7mm
\begin{picture}(5,5)(50,0)
\put(20,20){\line(1,-1){10}}
\put(40,20){\line(-1,-1){10}}
\put(40,20){\line(1,-1){10}}
\put(60,20){\line(-1,-1){10}}
\put(60,20){\line(1,-1){10}}

\qbezier(60, 20)(60,10)(70,10)

\put(63,18){$\searrow^1$}
\put(58,5){$\nwarrow^2$}

\put(20,20){\circle*{2}}
\put(40,20){\circle*{2}}
\put(60,20){\circle*{2}}

\put(30,10){\circle*{2}}
\put(50,10){\circle*{2}}
\put(70,10){\circle*{2}}
\end{picture}

\\

\begin{array}{cccccc}
 1 & \text{} & 1 & \text{} & 1 & \text{} \\
 \text{} & 3 & \text{} & 3 & \text{} & \frac{3}{2} \\
 2 & \text{} & 5 & \text{} & 5 & \text{} \\
 \text{} & 4 & \text{} & 7 & \text{} & \frac{7}{2} \\
 2 & \text{} & 6 & \text{} & 9 & \text{} \\
 \text{} & 4 & \text{} & 8 & \text{} & \frac{11}{2} \\
 2 & \text{} & 6 & \text{} & 10 & \text{} \\
 \text{} & 4 & \text{} & 8 & \text{} & \frac{9}{2} \\
 2 & \text{} & 6 & \text{} & 7 & \text{} \\
 \text{} & 4 & \text{} & 5 & \text{} & \frac{5}{2} \\
 2 & \text{} & 3 & \text{} & 3 & \text{} \\
 \text{} & 1 & \text{} & 1 & \text{} & \frac{1}{2} \\
\end{array}

\end{array}

&

C_6
\begin{array}{c}

\unitlength 0.7mm
\begin{picture}(5,5)(50,0)
\put(20,20){\line(1,-1){10}}
\put(40,20){\line(-1,-1){10}}
\put(40,20){\line(1,-1){10}}
\put(60,20){\line(-1,-1){10}}
\put(60,20){\line(1,-1){10}}

\qbezier(60, 20)(60,10)(70,10)

\put(63,18){$\searrow^2$}
\put(58,5){$\nwarrow^1$}

\put(20,20){\circle*{2}}
\put(40,20){\circle*{2}}
\put(60,20){\circle*{2}}

\put(30,10){\circle*{2}}
\put(50,10){\circle*{2}}
\put(70,10){\circle*{2}}
\end{picture}

\\

\begin{array}{cccccc}
 \frac{1}{2} & \text{ } & \frac{1}{2} & \text{ } & \frac{1}{2} & \text{ } \\
 \text{ } & \frac{3}{2} & \text{ } & \frac{3}{2} & \text{ } & 2 \\
 1 & \text{ } & \frac{5}{2} & \text{ } & 3 & \text{ } \\
 \text{ } & 2 & \text{ } & 4 & \text{ } & 4 \\
 1 & \text{ } & \frac{7}{2} & \text{ } & 5 & \text{ } \\
 \text{ } & \frac{5}{2} & \text{ } & \frac{9}{2} & \text{ } & 6 \\
 \frac{3}{2} & \text{ } & \frac{7}{2} & \text{ } & \frac{11}{2} & \text{ } \\
 \text{ } & \frac{5}{2} & \text{ } & \frac{9}{2} & \text{ } & 5 \\
 1 & \text{ } & \frac{7}{2} & \text{ } & 4 & \text{ } \\
 \text{ } & 2 & \text{ } & 3 & \text{ } & 3 \\
 1 & \text{ } & \frac{3}{2} & \text{ } & 2 & \text{ } \\
 \text{ } & \frac{1}{2} & \text{ } & \frac{1}{2} & \text{ } & 1 \\
\end{array}

\end{array}
\end{array}
$$
}
\caption{Values of the Weyl vector on the half-quiver of positive roots. Cases $A_6$, $D_6$, $B_6$, $C_6$.}
\label{tab:rhoA6B6C6D6} 
\end{table}

\bigskip

\subsection{Miscellaneous remarks}

\subsubsection*{Special case ${\mathcal G}={\mathcal A}_k(SU(2))$ and the path matrix.}

The path matrix\footnote{The terminology comes from the fact that its  elements give the dimensions of the vector spaces of essential paths \cite{Ocneanu:paths},\cite{Zuber:Bariloche}, on the fusion graph (also the dimensions of the vector spaces underlying the Gelfand-Ponomarev preprojective algebra associated with the corresponding unoriented quiver).} $X = \sum_{n=1}^r F_n$ can be defined for an arbitrary module-category associated with $SU(2)$ and specified by a simply-laced Dynkin diagram $G$.
The {\sl sum} over columns of $X$ gives the height vector of the Lie group defined by the same Dynkin diagram.
For $G=A_{r}$,  the Lie group is $SU(N)$, with $N=r+1=k+2$,  and $X$ is a $r\times r$ matrix of increasing concentric rectangular rings of integers starting from $1$ on the four sides (see the $A_{11}$ example given below). 
In that case we may  take all indices $a,b$ and $n$ from $1$ to $r$. The fusion coefficients $F_{p,m,n} = (F_p)_{m,n}$ are now completely symmetric with respect to permutation of their indices since the conjugation of representations is trivial for $SU(2)$.  Let then $\overset{\curvearrowleft}{X}$ be the matrix obtained from $X$ by cyclicly permuting the lines, \ie the matrix with entries 
$\overset{\curvearrowleft}{X}_{1,n} = X_{r,n}$ and $\overset{\curvearrowleft}{X}_{m,n} = X_{m-1,n}$ for $m>1$. 
From the previous symmetry property and using equation \eqref{bigprodeq} one obtains immediately, but for $G=A_{k+1}$ only,  the equality:
$$\prod_{m,n}  \;  (X+\overset{\curvearrowleft}{X})_{m,n} \, = 2^r \, \mathrm{\sfac}(r)^2$$
We display the matrix $F_4$ at level $10$ and the path matrix $X$; the Dynkin diagram is $A_{11}$.
{\tiny
$$
F_4=
\left(
\begin{array}{ccccccccccc}
 0 & 0 & 0 & 1 & 0 & 0 & 0 & 0 & 0 & 0 & 0 \\
 0 & 0 & 1 & 0 & 1 & 0 & 0 & 0 & 0 & 0 & 0 \\
 0 & 1 & 0 & 1 & 0 & 1 & 0 & 0 & 0 & 0 & 0 \\
 1 & 0 & 1 & 0 & 1 & 0 & 1 & 0 & 0 & 0 & 0 \\
 0 & 1 & 0 & 1 & 0 & 1 & 0 & 1 & 0 & 0 & 0 \\
 0 & 0 & 1 & 0 & 1 & 0 & 1 & 0 & 1 & 0 & 0 \\
 0 & 0 & 0 & 1 & 0 & 1 & 0 & 1 & 0 & 1 & 0 \\
 0 & 0 & 0 & 0 & 1 & 0 & 1 & 0 & 1 & 0 & 1 \\
 0 & 0 & 0 & 0 & 0 & 1 & 0 & 1 & 0 & 1 & 0 \\
 0 & 0 & 0 & 0 & 0 & 0 & 1 & 0 & 1 & 0 & 0 \\
 0 & 0 & 0 & 0 & 0 & 0 & 0 & 1 & 0 & 0 & 0 \\
\end{array}
\right)
\quad
X=
\left(
\begin{array}{ccccccccccc}
 1 & 1 & 1 & 1 & 1 & 1 & 1 & 1 & 1 & 1 & 1 \\
 1 & 2 & 2 & 2 & 2 & 2 & 2 & 2 & 2 & 2 & 1 \\
 1 & 2 & 3 & 3 & 3 & 3 & 3 & 3 & 3 & 2 & 1 \\
 1 & 2 & 3 & 4 & 4 & 4 & 4 & 4 & 3 & 2 & 1 \\
 1 & 2 & 3 & 4 & 5 & 5 & 5 & 4 & 3 & 2 & 1 \\
 1 & 2 & 3 & 4 & 5 & 6 & 5 & 4 & 3 & 2 & 1 \\
 1 & 2 & 3 & 4 & 5 & 5 & 5 & 4 & 3 & 2 & 1 \\
 1 & 2 & 3 & 4 & 4 & 4 & 4 & 4 & 3 & 2 & 1 \\
 1 & 2 & 3 & 3 & 3 & 3 & 3 & 3 & 3 & 2 & 1 \\
 1 & 2 & 2 & 2 & 2 & 2 & 2 & 2 & 2 & 2 & 1 \\
 1 & 1 & 1 & 1 & 1 & 1 & 1 & 1 & 1 & 1 & 1 \\
\end{array}
\right)
$$
}

\subsubsection*{Fourier-like considerations.}

Relating the Weyl denominator of a simple Lie group to fusion matrices $F_n$ describing module-categories of type $SU(2)$ at level $k$, leads, if we diagonalize those matrices, to non-trivial identities relating  trigonometric lines, sums of Chebyshev matrix polynomials, and superfactorial functions.
Using $F_n = \Psi^\dag \, \text{diag}(\Delta_n) \, \Psi$  -- see section \ref{sec:fusiontoscalarproducts} --  and equation \eqref{bigprodeq}, one finds immediately
$$\prod_{s=1}^{N} \prod_{a=1}^r \sum_{b=1}^r  \; \left(\Psi^\dag \; \text{diag}\left[\frac{\sin \left(\frac{\pi \epsilon_j}{2 N }+\frac{\pi  s
  \epsilon_j}{N }\right)}{\sin \left(\frac{\pi \epsilon_j}{2 N
   }\right)}\right]_{\epsilon_j \in {\mathcal E}xp} \; \Psi\right)_{a,b} \; = \;  (\mathrm{\sfac}_G)^2$$ 
  
\noindent
 For instance, with $G=E_6$,  ${\mathcal E}xp={1,4,5,7,8,11}$, $N=12$,  and $\Psi$ as below, one can check that the above expression is indeed equal to $(1! \, 4! \, 5! \,  7! \,  8!  \, 11!)^2$.
{\tiny $$ \Psi= \frac{1}{2 \sqrt{2}\, } \; \left(
\begin{array}{cccccc}
\phi_{-} &\phi_{+} & \sqrt{2}\, \phi_{+} &\phi_{+} &\phi_{-} &
   \sqrt{2}\, \phi_{-} \\
 \sqrt{2}\,  & \sqrt{2}\,  & 0 & -\sqrt{2}\,  & -\sqrt{2}\,  & 0 \\
\phi_{+} &\phi_{-} & -\sqrt{2}\, \phi_{-} &\phi_{-} &\phi_{+} &
   -\sqrt{2}\, \phi_{+} \\
\phi_{+} & -\phi_{-} & -\sqrt{2}\, \phi_{-} & -\phi_{-} &\phi_{+} &
   \sqrt{2}\, \phi_{+} \\
 \sqrt{2}\,  & -\sqrt{2}\,  & 0 & \sqrt{2}\,  & -\sqrt{2}\,  & 0 \\
\phi_{-} & -\phi_{+} & \sqrt{2}\, \phi_{+} & -\phi_{+} &\phi_{-} &
   -\sqrt{2}\, \phi_{-} \\
\end{array}
\right) \quad \text{with} \quad \phi_{\pm}=\sqrt{\frac{1}{3} \left(3\pm\sqrt{3}\right)}$$}

One could think of using equation\eqref{bigprodeq} for an arbitrary Dynkin diagram $G$ in order to show that the Weyl denominator can indeed be expressed in terms of  a Lie superfactorial function,  
but the direct evaluation of the left hand side, even for $A_r$, with $r$ generic, does not look straightforward.
In the $A_r$ cases one can take for $\Psi$ the modular $S$-matrix with entries 
$S_{m,n}=  \sqrt{\frac{2}{N }} \sin \left(\frac{\pi \, m \, n}{N }\right)$
where $N=r+1$.  Equation\eqref{bigprodeq} then reads:
$$
\prod_{s, m = 1}^{N -1}\,
\sum _{n, j =1}^{N -1}
\; \frac{2}{N} \;  
\sin \left(\frac{\pi  j m}{N }\right) 
\sin \left(\frac{\pi  j n}{N }\right) 
 \dfrac
 {\sin \left(\dfrac{\pi  j (2 s+1)}{2N }\right)}
 {\sin \left(\dfrac{\pi  j}{2 N }\right) }  
 \; = \; (\mathrm{\sfac}[N-1])^2 = \prod_{s=1}^{N-1} \, s!^2 $$

\subsubsection*{Evaluation of the fusion product formula  for $A_r$}
\label{sec:fusionmethodformulaforA}

Another possibility, for generic $r$,  is to use the last identity of section \ref{sec:weyldenominatorfromfusion}{} and directly evaluate its right hand side. We consider the case $A_r$.  The adjacency matrix $G=F_2$ has matrix elements $G_{ab}=\delta_{a+1,b} + \delta_{a-1,b}$ where $\delta$ is the Kronecker symbol.
For $b=1\ldots r$ and $n=1\ldots N=r+1$, call $\lambda_{n,b}=\sum_{a=1}^r(F_{n}+F{_{n-1}})_{a,b}$.
Assuming $N=r+1$ even, one shows immediately that, if $n\leq N/2$, $\lambda_{n,b} = 2b$ for $b=1,2,\ldots, n-1$,
$\lambda_{n,b} = 2n-1$ for $b=n,n+1,\ldots, N-n$,
and $\lambda_{n,b} = 2(N-b)$ for $b=N-n+1,\ldots, N-1$.
Moreover we have the symmetry properties $\lambda_{n,b}=\lambda_{N-n+1,b}$ and $\lambda_{n,b}=\lambda_{n,N-b}$.
If  $n\leq N/2$, we have  $\prod_{b=1}^{N-1} \lambda_{n,b} = (\prod_{b=1}^{n-1} (2b))^2 \times (2n-1)^{N-2n+1} = 4^{n-1} (n-1)!^2  (2n-1)^{N-2n+1}$.
Because of the symmetry properties one obtain the same value if $n\geq N/2+1$.
The last step is the product over $n$:
${\mathcal D}^2 = \prod_{n=1}^{N} \prod_{b=1}^{N-1} \lambda_{n,b} =  (\prod_{n=1}^{N/2} \prod_{b=1}^{N-1} \lambda_{n,b})^2$, so 
${\mathcal D} =  \prod_{n=1}^{N/2}\prod_{b=1}^{N-1} \lambda_{n,b}= {\mathcal D}_1 {\mathcal D}_2 {\mathcal D}_3$, with  
${\mathcal D}_1=\prod_{n=1}^{N/2}  4^{n-1}$, ${\mathcal D}_2=  \prod_{n=1}^{N/2}  (n-1)!^2$,
${\mathcal D}_3= \prod_{n=1}^{N/2}  (2n-1)^{N-(2n-1)}$.
Evaluation of the first two brackets is trivial:  ${\mathcal D}_1= 2^{\frac{N}{2}(\frac{N}{2}-1)}$,  ${\mathcal D}_2=\mathrm{\sfac}(\frac{N}{2} - 1)^2$. 
Evaluation of the third  being more involved, we give some details: 
let us write ${\mathcal D}_3$ as $u/v$, with $u =( \prod_{n=1}^{N/2} (2n-1))^N$ and $v= \prod_{n=1}^{N/2} (2n-1)^{2n-1}$. 
The numerator $u$ is the $N$-th power of a double factorial that we write, $N=2s$ being even,  
$u=(2^{s} \;  \Gamma(\frac{2s+1}{2}) / \sqrt \pi)^{2s}$, where we introduced the Euler Gamma function.
We can express the denominator $v$ in terms of the hyperfactorial function $\text{H}(n)=\prod_{k=1}^n k^k$, see \cite{SloanePlouffe}, 
and obtain $v=2^{s^2} \pi^{-1/2} \; \text{G}(3/2)^2 \; \text{H}(\frac{2s-1}{2})^2$, 
where $\text{G}$ is the Barnes $\text{G}$-function.
Using the identity $\text{H}(z)= \Gamma(z+1)^z/\text{G}(z+1)$ for $z={(2s-1)}/{2}$
 and the known\footnote{As the result for ${\mathcal D}$  was already obtained by a different method, those identities found for $\text{G}((2s+1)/2)$ constitute a consistency check} values of $\text{G}((2s+1)/2)$, see \cite{Weisstein:BarnesG}, that we prefer to write 
$$\text{G}(\frac{2s+1}{2})=\frac{\text{G}\left(\frac{3}{2}\right)  \pi ^{\tfrac{s}{2}} \sqrt{\mathrm{\sfac}(2s-1)}}{2^{s(s-1)} \sqrt{\pi}  \; \mathrm{\sfac}(s-1)  \; \sqrt{(2 s-1)\text{!!}}}$$
with $\mathrm{\sfac}(z)=\text{G}(z+2)$ if $z$ is a positive integer, one obtains
${\mathcal D}_3= 2^{\frac{-N}{2}(\frac{N}{2}-1)} \, \mathrm{\sfac}(N-1)/ \mathrm{\sfac}(\frac{N}{2}-1)^2$, hence the expected result: ${\mathcal D}=\mathrm{\sfac}(N-1)$.

\noindent
Details differ if we assume that $N$ is odd but the calculation is essentially the same.

\smallskip
For Dynkin diagrams of type $B_r$, $C_r$, $D_r$, one may proceed in a similar manner, starting from their adjacency matrices.
Notice that, from the very beginning, one has the choice of performing the products, first over $b$, then over $n$, or in the opposite order (this latter possibility  amounts to start the calculation by taking the product over the orbits of the Coxeter element), details then differ. In all cases the calculation is as straightforward, but as cumbersome as for $A_r$, and the result for ${\mathcal D}$ is known anyway (cf. Appendix). We leave this as an exercise.
\section*{Appendix. {\small From the Weyl denominator to superfactorials (classical series)}}
\label{sec:genericproofforABCD}

In the  $A_r$ case,  the familiar Schur  formula used to calculate dimensions of irreducible representations described by Young tableaux involves a division by a product of successive factorials. Interpretation in terms of exponents of the group $SU(r+1)$ is then immediate since the latter are given by the consecutive integers ${1,2,\ldots r}$. 
For the other classical series $B_r, C_r, D_r$,  the fact that one can express  the (classical) Weyl denominator ${\mathcal D}$ as a product of factorials seems also to belong to the folklore, nevertheless we did not find any precise reference -- except  \cite{CordobaRC} -- relating this property to exponents of Lie groups.
 It was in the spirit of the paper to obtain the scalar products between roots and the Weyl vector from the fusion coefficients of appropriate module-categories of type $SU(2)$. 
 This method, which is both fast and straightforward, can be used for every particular case, in particular for exceptional Lie groups and it was described in the main body of the paper.
However, for classical series with $r$ generic,  we shall proceed differently in this appendix and present a proof that has the advantage of being quite elementary, in the sense that it only uses well-known Lie group theoretical results.
In the approach used below we  determine explicitly the list of values of the needed scalar products simply by
expressing both fundamental weights and positive roots in terms of the elements of the canonical basis of a real euclidean space, and performing the actual calculation.

\bigskip 
For $B_r, C_r, D_r$, using well-known expressions, we write both fundamental weights and positive roots in terms of the elements $\epsilon_i$ of the canonical basis of $\RR^r$.
In the case $A_r$, it is  convenient to embed the root space in $\RR^{r+1}$, or equivalently in the space of diagonal matrices $(r+1)\times(r+1)$. 
The calculations are straightforward and only use the fact that $\langle \epsilon_i, \epsilon_j \rangle = \delta_{i,j}$, where $\delta_{i,j}$ is the Kronecker symbol.

\smallskip
\noindent
We introduce the following notation:

$[a, b]  =  (a, a+1, a+2, \ldots, b-1, b) \quad \text{if} \quad a \leq b \quad \text{or} \quad  (a, a-1, a-2, \ldots, b+1, b) \quad \text{if} \quad a \geq b$
Here, $a$ and $b$ are positive rationals with $b-a$ an integer. We also set $[a,a]=(a)$.

\bigskip
\noindent
{\bf Case $A_r$.} 

\noindent
Fundamental weights:
$\omega_i= \epsilon_1+ \ldots  + \epsilon_i-\frac{i}{r+1}\,  \sigma \quad \text{where} \quad \sigma = \sum_{i=1}^{r+1} \, \epsilon_i$

\noindent
Weyl vector:
$\rho = \sum_{i=1}^{r} \omega_i = \Sigma_{p=0}^{r-1}\, (r-p)  \epsilon_{p+1} - \frac{r}{2} \, \sigma$

\noindent
Positive roots:
$\alpha_{i, j} = \epsilon_i - \epsilon_j \, : 1 \leq i < j \leq r+1$

\noindent
Scalar products:
$\langle \rho,  \alpha_{i, j} \rangle = \sum_{p=1}^r (\frac{r}{2} - (p-1)) (\delta_{p,i} - \delta_{p,j}) - \frac{r}{2} (\delta_{r+1,i} - \delta_{r+1,j})$

\noindent
Hence, we obtain the following tuples:
$$\text{for} \, i=1,\ldots r, \quad \left( \langle \rho,  \alpha_{i, j} \rangle\right)_{j= i+1, \ldots, r+1} = [1,r-(i-1)]$$

\noindent
Therefore ${\mathcal D} = \Pi_{i=1}^r r! =  \mathrm{\sfac}_{A_r} [1] = \mathrm{\sfac} \, (r)$

\bigskip
\noindent
{\bf Case $B_r$.} 

\noindent
Fundamental weights:
$\omega_i= \epsilon_1 + \ldots + \epsilon_i, \quad i \neq r \quad \text{and} \quad \omega_r = \frac{1}{2} \, (\epsilon_1 + \ldots + \epsilon_r)$

\noindent
Weyl vector:
$\rho = \sum_{i=1}^{r} \omega_i = \Sigma_{p=1}^{r-1}\, (r-p)  \epsilon_{p} + \frac{r}{2} \, (\epsilon_1 + \ldots + \epsilon_r)$

\noindent
Positive roots\footnote{Those called $\alpha_{i}^{(3)}$ are short, the others are long.} :\\
$\alpha_{i, j}^{(1)} = \epsilon_i + \epsilon_j \, : 1 \leq i < j \leq r$, $\alpha_{i, j}^{(2)} = \epsilon_i - \epsilon_j \, : 1 \leq i < j \leq r$ and $\alpha_{i}^{(3)} = \epsilon_i  : 1 \leq i  \leq r$

\noindent
Scalar products:\\
$\langle \rho,  \alpha_{i, j}^{(1)} \rangle = \sum_{p=1}^{r-1} (r - p) (\delta_{p,i} + \delta_{p,j}) + \frac{1}{2} \sum_{s=1}^r (\delta_{i,s} + \delta_{j,s})$ \\
$\langle \rho,  \alpha_{i, j}^{(2)} \rangle = \sum_{p=1}^{r-1} (r - p) (\delta_{p,i} - \delta_{p,j}) + \frac{1}{2} \sum_{s=1}^r (\delta_{i,s} - \delta_{j,s})$ \\
$\langle \rho,  \alpha_{i}^{(3)} \rangle = \sum_{p=1}^{r-1} (r - p)  \delta_{p,i} + \frac{1}{2}  \sum_{s=1}^r \delta_{s,i}$

\noindent
Hence, we obtain the following tuples:
\begin{eqnarray*}
\text{for} \, i=1,\ldots r, \quad \left( \langle \rho,  \alpha_{i, j}^{(1)} \rangle\right)_{j= i+1, \ldots, r+1} &=& [2r-2, r], [2r-4, r-1], \ldots, (6,5,4), (4,3), (2)\\
\text{for} \, i=1,\ldots r, \quad \left( \langle \rho,  \alpha_{i, j}^{(2)} \rangle\right)_{j= i+1, \ldots, r+1} &=& [1,r-1], [1, r-2], [1, r-3], \ldots, (1,2,3), (1,2), (1)\\
\text{for} \, i=1,\ldots r, \quad \left( \langle \rho,  \alpha_{i}^{(3)} \rangle\right)_{j= i+1, \ldots, r+1} &=& ((2r-1)/2, (2r-3)/2, \ldots, 3/2, 1/2) =  [(2r-1)/2, 1/2]
\end{eqnarray*}

Collecting together and multiplying the values obtained for the tuples $\alpha_{i, j}^{(1)}$ and $\alpha_{i, j}^{(2)}$ with the same index $i$, one obtains a product of factorials $s!$ for $s= 2r-2, 2r-4, \ldots, 6, 4, 2$. Together with the contribution of  $\alpha_{i, j}^{(3)}$ one obtains  ${\mathcal D} = \dfrac{1}{2^r} \,  \Pi_{s \in {\mathcal E}}\,  {s!}$ where $s = 1, 3, 5, \ldots, 2r-1$ are the exponents of $B_r$. 
However,  it is better to write this result while keeping the individual values of the scalar products since the quantum Weyl denominator is obtained from the classical one by replacing those individual values by the corresponding $q$-numbers. So we prefer to write:
${\mathcal D} = \Pi_{s \in {\mathcal E}}\,  \widetilde{s!}$ with $\widetilde {s!} = (s/2)(s-1)(s-2)\ldots (2)(1)$. Therefore  ${\mathcal D} = \mathrm{\sfac}_{B_r} [1]$.

\bigskip
\noindent
{\bf Case $C_r$.} 

\noindent
The  factor $\sqrt 2$ appearing below in the given expressions for fundamental weights $\omega_i$ and positive roots $\alpha_{i,j}$ is introduced to ensure that the norm square of long roots is $2$.

\noindent
Fundamental weights:
$\sqrt 2 \; \omega_i= \epsilon_1 + \ldots + \epsilon_i$

\noindent
Weyl vector:
$\rho = \sum_{i=1}^{r} \omega_i = (1/\sqrt 2) \; \Sigma_{p=1}^{r}\, (r-(p-1))  \epsilon_{p}$

\noindent
Positive roots\footnote{The long roots $\alpha_{i}^{(3)}$ have norm square $2$. The others are short (their norm square is $1$).} 
$\sqrt{2}\, \alpha_{i, j}^{(1)} = (\epsilon_i + \epsilon_j) \, : 1 \leq i < j \leq r$; $\sqrt{2}\, \alpha_{i, j}^{(2)} = (\epsilon_i - \epsilon_j) \, : 1 \leq i < j \leq r$; $\sqrt{2}\, \alpha_{i}^{(3)} = 2 \epsilon_i  : 1 \leq i  \leq r$

\noindent
Scalar products:\\
$\langle \rho,  \alpha_{i, j}^{(1)} \rangle = (1/2) \sum_{p=1}^{r} (r - (p-1)) (\delta_{p,i} + \delta_{p,j})$ \\
$\langle \rho,  \alpha_{i, j}^{(2)} \rangle = (1/2) \sum_{p=1}^{r} (r - (p-1)) (\delta_{p,i} - \delta_{p,j})$ \\
$\langle \rho,  \alpha_{i}^{(3)} \rangle = (1/2) \sum_{p=1}^{r} (r - (p-1)) \,  2 \delta_{p,i}$

\noindent
Hence, we obtain the following tuples:
\begin{eqnarray*}
\text{for} \, i=1,\ldots r, \quad \left( \langle \rho,  \alpha_{i, j}^{(1)} \rangle\right)_{j= i+1, \ldots, r+1} &=& [(2r-1)/2, (r+1)/2], [(2r-3)/2, r/2],  \ldots, (5/2,4/2), (3/2)\\
\text{for} \, i=1,\ldots r, \quad \left( \langle \rho,  \alpha_{i, j}^{(2)} \rangle\right)_{j= i+1, \ldots, r+1} &=& [1/2,(r-1)/2], [1/2, (r-2)/2],\ldots, (1/2,2/2), (1/2)\\
\text{for} \, i=1,\ldots r, \quad  \langle \rho,  \alpha_{i}^{(3)} \rangle  &=& (2r/2, (2r-2)/2, \ldots, 4/2, 2/2) =  [r, 1]
\end{eqnarray*}

Collecting the values obtained for the tuples $\alpha_{i, j}^{(1)}$, $\alpha_{i, j}^{(2)}$ and $\alpha_{i, j}^{(3)}$ with the same index $i=k-1$, one obtains
tuples $[(2r-(2k+1))/2, (r-(k-1))/2], (r-k), [(r-(k+1))/2, 1/2]$, for $k=0, \ldots, r-2$. Multiplying their entries together for a given $k$, one obtains $\dfrac{1}{2^{2(r-(1+k))}} (2r-(2k+1))!$.  Hence 
${\mathcal D} = \frac{1}{2^{r(r-1)}} \Pi_{s \in {\mathcal E}}\,  {s!}$ where $s = 1, 3, 5, \ldots, 2r-1$ are the exponents of $C_r$.
For the same reason as before, we prefer to write:\\
${\mathcal D} = \Pi_{s \in {\mathcal E}}\,  \widetilde{s!}$ with $\widetilde {s!} = \dfrac{s}{2} \dfrac{s-1}{2} \ldots \dfrac{s-(s-3)/2}{2}  \; (s- \dfrac{s-1}{2}) \; \dfrac{s-(s+1)/2}{2} \ldots \dfrac{2}{2} \dfrac{1}{2}$.
\\
Therefore  ${\mathcal D} = \mathrm{\sfac}_{C_r} [1]$.   

\vskip 1.cm
\bigskip
\noindent
{\bf Case $D_r$.} 

\noindent
Fundamental weights: 
 $\omega_i= \epsilon_1 + \ldots + \epsilon_i$ for $i=1,\ldots r-2$. \\ The last two are 
$\omega_{r-1}=\dfrac{1}{2} ( \epsilon_1 + \ldots + \epsilon_{r-1} - \epsilon_r)$ and
$\omega_{r}= \dfrac{1}{2}  ( \epsilon_1 + \ldots + \epsilon_{r-1} + \epsilon_r)$.

\noindent
Weyl vector:
$\rho = \sum_{i=1}^{r} \omega_i = \Sigma_{p=1}^{r-1}\, (r-p)  \epsilon_{p}$

\noindent
Positive roots:
$\alpha_{i, j}^{(1)} = \epsilon_i + \epsilon_j \, : 1 \leq i < j \leq r$ and $\alpha_{i, j}^{(2)} = \epsilon_i - \epsilon_j \, : 1 \leq i < j \leq r$

\noindent
Scalar products:
$\langle \rho,  \alpha_{i, j}^{(1)} \rangle = \sum_{p=1}^{r-1} (r - p) (\delta_{p,i} + \delta_{p,j})$, and
$\langle \rho,  \alpha_{i, j}^{(2)} \rangle = \sum_{p=1}^{r-1} (r - p) (\delta_{p,i} - \delta_{p,j}) $

\noindent
Hence, we obtain the following tuples:
\begin{eqnarray*}
\text{for} \, i=1,\ldots r, \quad \left( \langle \rho,  \alpha_{i, j}^{(1)} \rangle\right)_{j= i+1, \ldots, r} &=& [2r-3,r-1], [2r-5,r-2], \ldots, (5,4,3), (3,2), (1)\\
\text{for} \, i=1,\ldots r, \quad \left( \langle \rho,  \alpha_{i, j}^{(2)} \rangle\right)_{j= i+1, \ldots, r} &=& [1, r-1], [1, r-2], \ldots, (1,2,3), (1,2), (1)\\
\end{eqnarray*}

Collecting the values obtained for the tuples $\alpha_{i, j}^{(1)}$ and $\alpha_{i+1, j}^{(2)}$, for a given $i$, with $i=1, \ldots, r$,  one builds a tuple $[2r-(2i+1), 1]$.
Multiplying their entries together one obtains a contribution $(2r-(2i+1))!$ to the Weyl denominator.  Finally the values $\alpha_{i = 1, j}^{(2)} = [1, r-1]$ give an extra factor $(r-1)!$.
Hence ${\mathcal D} = (r-1)!  \; \Pi_{s=1}^{r-1} \, (2s-1)! =  \Pi_{s \in {\mathcal E}}\,  {s!}$ where $s = 1, 3, 5, \ldots, 2r-3; r-1$ are the exponents of $D_r$ (when $r$ is even, the exponent $r-1$ appears with multiplicity $2$ in ${\mathcal E}$).  Therefore  ${\mathcal D} = \mathrm{\sfac}_{D_r} [1]$. 

\bigskip
For the exceptional cases, the expressions giving fundamental weights and positive roots in terms of an orthonormal basis are more involved, but they are known (see for instance \cite{ConwaySloane}). We could use the same technique as above to obtain the needed lists of scalar products that determine both ${\mathcal D}$ and ${\mathcal D}_q$.
We nevertheless leave this as an exercise since the lists have been  obtained by another technique (correspondence with fusion coefficients) in section \ref{sec:periodicquiver}.

 \end{document}